\documentclass{article}
\usepackage[english]{babel}
\usepackage{amsmath,amssymb,stmaryrd,enumerate,bbm,latexsym,theorem}
\usepackage[utf8]{inputenc}
\usepackage{graphicx}
\usepackage{float}

\usepackage[top=2.5cm,bottom=2.5cm,right=2cm,left=2cm]{geometry}
\usepackage{authblk}
\numberwithin{equation}{section}

\newcommand{\assign}{:=}

\newcommand{\tmmathbf}[1]{\ensuremath{\boldsymbol{#1}}}

\newcommand{\tmop}[1]{\ensuremath{\operatorname{#1}}}

\newcommand{\tmtextit}[1]{{\itshape{#1}}}
\newenvironment{enumeratenumeric}{\begin{enumerate}[1.] }{\end{enumerate}}
\newenvironment{itemizedot}{\begin{itemize} }{\end{itemize}}
\newenvironment{proof}{\noindent\textbf{Proof\ }}{\hspace*{\fill}$\Box$\medskip}
\newcounter{nnacknowledgments}

{\theorembodyfont{\rmfamily}\newtheorem{acknowledgments*}[nnacknowledgments]{Acknowledgments}}
\newtheorem{theorem}{Theorem}[section]

\newtheorem{proposition}[theorem]{Proposition}

\usepackage[hidelinks]{hyperref}

\date{}
\author{Tej-Eddine Ghoul \thanks{email: teg6@nyu.edu}}
\author{Nader Masmoudi \thanks{email: nm30@nyu.edu}}
\author{Eliot Pacherie  \thanks{email: ep2699@nyu.edu}}
\affil{NYUAD Research Institute, New York University Abu Dhabi}

\begin{document}

\title{Nonlinear enhanced dissipation in viscous Burgers type equations II}

\maketitle

\begin{abstract}
  In this follow up but self contained paper, we focus on the viscous Burgers equation. There,
  using the Hopf-Cole transformation, we compute the long time behavior of
  solutions for some classes of infinite mass initial datas. We show that an
  enhanced dissipation effect occurs generically, that is the decay rate in
  time is better than if we considered instead the heat equations for the same
  inital value. We also show the existence of a kind of global attractor per
  class. 
\end{abstract}

\section{Introduction and presentation of the results}

We are interested in this paper on the long time behavior of solutions to
either the heat or the viscous Burgers equation on the real line. It is well
known that for the heat equation $\partial_t u - \partial_x^2 u = 0$, for an
initial data $u_0 \in L^1 (\mathbb{R})$ we have the asymptotic profile
\[ \sqrt{t} u \left( z \sqrt{t}, t \right) \rightarrow \frac{\int_{\mathbb{R}}
   u_0}{\sqrt{4 \pi}} e^{- \frac{z^2}{4}} \]
when $t \rightarrow + \infty$, uniformly in $z \in \mathbb{R}$.

A similar result holds for the viscous Burgers equation $\partial_t v -
\partial_x^2 v + v \partial_x v = 0$ for initial data $v_0 \in L^1
(\mathbb{R})$, (see {\cite{MR1124296}}, {\cite{MR3583525}},
{\cite{MR2313029}}), as we have
\[ \sqrt{t} v \left( z \sqrt{t}, t \right) \rightarrow \frac{2 \left(
   e^{\frac{M}{2}} - 1 \right) e^{- z^2 / 4}}{e^{\frac{M}{2}} \sqrt{4 \pi} +
   \left( 1 - e^{\frac{M}{2}} \right) \int_{- \infty}^z e^{- s^2 / 4} d s} \]
when $t \rightarrow + \infty$, uniformly in $z \in \mathbb{R}$, with $M =
\int_{\mathbb{R}} v_0$.

In both case, the $L^{\infty}$ norm of the solution decays like $t^{-
\frac{1}{2}}$. Other asymptotic behavior results have been established in
other convection-diffusion equations for initial datas in $L^1 (\mathbb{R})$,
we refer to {\cite{arxivzuazuabis}} and references therein, as well as
{\cite{MR2666480}}, {\cite{MR2172561}}.

\

Our first goal is to show that, for some class of initial datas that are not
in $L^1 (\mathbb{R})$, the decay rate for the viscous Burgers equation is
better than the one of the heat equation.

\subsection{Nonlinear enhanced dissipation}

Enhanced dissipation results are well known for the heat equation with an
additional linear transport term (see for instance {\cite{MR3621820}},
{\cite{MR2434887}}, {\cite{MR4156602}}, and references therein) or for
Navier-Stokes on $\mathbb{T} \times \mathbb{R}$ (see {\cite{MR3448924}},
{\cite{MR4134940}}, {\cite{MR4030287}}, {\cite{MR4176913}}). It has also been
established in {\cite{MR1266100}} for some convection-diffusion equations on
the real line.

\

The main result of this subsection is an estimation of the decay in time for
solutions to the viscous Burgers equation in a particular class of initial
datas, on which the decay is better than the one for the heat equation.

\begin{theorem}
  \label{mainth}Consider the problem
  \begin{equation}
    \partial_t f - \partial_x^2 f + f \partial_x f = 0 \label{vb}
  \end{equation}
  for an initial data $f_0 \in C^0 (\mathbb{R}, \mathbb{R})$, and suppose that
  there exists $\kappa_1, \kappa_2 > 0$ and $\alpha \in] 0, 1 [$ such that
  \[  \frac{\kappa_2}{(1 + | y |)^{\alpha}} \geqslant f_0 (y) \geqslant
     \frac{\kappa_1}{(1 + | y |)^{\alpha}} . \]
  Then, there exists $K_1, K_2 > 0$ depending on $\kappa_1, \kappa_2, \alpha,$
  such that, for all $t \geqslant 0$,
  \[ \frac{K_1}{(1 + t)^{\frac{\alpha}{1 + \alpha}}} \leqslant \| f (., t)
     \|_{L^{\infty} (\mathbb{R})} \leqslant \frac{K_2}{(1 +
     t)^{\frac{\alpha}{1 + \alpha}}} . \]
\end{theorem}

This result is to be compared with the decay estimates for solutions to the
heat equation for the same initial conditions.

\begin{proposition}
  \label{prop12}Consider the problem
  \[ \partial_t f - \partial_x^2 f = 0 \]
  for an initial data $f_0 \in C^0 (\mathbb{R}, \mathbb{R})$ satisfying the
  conditions of Theorem \ref{mainth}. Then, there exists $K_1, K_2 > 0$
  depending on $\kappa_1, \kappa_2, \alpha$ such that, for $t \geqslant 0$,
  \[ \frac{K_1}{(1 + t)^{\frac{\alpha}{2}}} \leqslant \| f (., t)
     \|_{L^{\infty} (\mathbb{R})} \leqslant \frac{K_2}{(1 +
     t)^{\frac{\alpha}{2}}} . \]
\end{proposition}

We see here that the dissipation is enhanced by adding the Burgers term $f
\partial_x f$ ($\frac{\alpha}{1 + \alpha} > \frac{\alpha}{2}$ for $\alpha \in]
0, 1 [$). It was proven in {\cite{GMP2}} that for a generalized version of the
viscous Burgers equation, this effect appears for a particular function (and
small perturbations of it) behaving like $\kappa_{\pm} | x |^{- \alpha}$ when
$x \rightarrow \pm \infty$. Theorem \ref{mainth} shows that this is in fact
true for all functions with this behavior at infinity for the viscous Burgers
equations.

If we suppose some decay on derivatives of the initial data, it is possible to
estimate the derivatives of the solution of the viscous Burgers equation
(\ref{vb}).

\begin{proposition}
  \label{Prop14}Consider the problem
  \[ \partial_t f - \partial_x^2 f + f \partial_x f = 0 \]
  for an initial data $f_0 \in C^j (\mathbb{R}, \mathbb{R})$ with $j \geqslant
  1$, and suppose that there exists $\kappa_1, \kappa_2 > 0$ and $\alpha \in]
  0, 1 [$ such that
  \[  \frac{\kappa_2}{(1 + | y |)^{\alpha}} \geqslant f_0 (y) \geqslant
     \frac{\kappa_1}{(1 + | y |)^{\alpha}} . \]
  Furthermore, if for any $i \in \mathbb{N}, i \leqslant j,$ there exists
  $\lambda_i > 0$ such that
  \[ | f_0^{(i)} (y) | \leqslant \frac{\lambda_i}{(1 + | y |)^{\alpha + i}},
  \]
  then for any $n, k$ such that $2 n + k \leqslant j$, there exists $K_{n, k}
  > 0$ depending on $\kappa_1, \kappa_2, (\lambda_i), \alpha, n$ and $k$ such
  that, for $t \geqslant 0$,
  \[ \| \partial_t^n \partial_x^k f (., t) \|_{L^{\infty} (\mathbb{R})}
     \leqslant \frac{K_{n, k}}{(1 + t)^{\frac{\alpha}{1 + \alpha} (1 + 2 n +
     k)}} . \]
\end{proposition}

This also has to be compared with the similar result for the heat equation,
where derivatives adds more decay than for viscous Burgers.

\begin{proposition}
  \label{prop15}Consider the problem
  \[ \partial_t f - \partial_x^2 f = 0 \]
  for an initial data $f_0 \in C^j (\mathbb{R}, \mathbb{R}), j \geqslant 1$
  satisfying the conditions of Proposition \ref{Prop14}. Then for any $n, k
  \in \mathbb{N}$, there exists $K_{n, k} > 0$ depending on $\kappa_1,
  \kappa_2, \alpha, n, k$ such that, for $t \geqslant 0, n + k \leqslant j$,
  \[ \| \partial_t^n \partial_x^k f (., t) \|_{L^{\infty} (\mathbb{R})}
     \leqslant \frac{K_{n, k}}{(1 + t)^{\frac{\alpha}{2} + n + \frac{k}{2}}} .
  \]
\end{proposition}

Remark that for the heat equation, adding a derivative in time or in position
add respectively an additional decay like $t^{- 1}$ and $t^{- \frac{1}{2}}$,
but in the viscous Burgers equation (\ref{vb}), these decays are respectively
$t^{- \frac{2 \alpha}{1 + \alpha}}$ and $t^{- \frac{\alpha}{1 + \alpha}}$. It
is possible that these rates are not optimal. See section \ref{54ss13} for the
strategy of the proof.

\subsection{Computation of the asymptotic profile}\label{s111}

We are interested in computing the limit when $t \rightarrow + \infty$ of
$t^{\frac{\alpha}{1 + \alpha}} f \left( z t^{\frac{1}{1 + \alpha}}, t \right)$
for some fixed $z \in \mathbb{R}$ where $f$ is the solution of the viscous
Burgers equation (\ref{vb}) with initial data $f_0$ similar to the one of
Theorem \ref{mainth}. In {\cite{GMP2}}, we computed the asymptotic profile in
one particular case, and we will see here that it is in fact, in some sense, a
global attractor. This will involve two particular functions defined by an
implicit equation. Consider the function
\[ g (y) \assign y + \frac{\kappa}{| y |^{\alpha}} . \]
We are interested in solving the equation $z = g (y (z))$ for the unknown
function $y$. Here is the plot of $\mathcal{C} \assign \{ (z, y) \in
\mathbb{R}^2, z = g (y) \}$ (in the case $\kappa = 1, \alpha = 1 / 3$):

\begin{figure}[H]
    \centering
    \includegraphics[width=20cm]{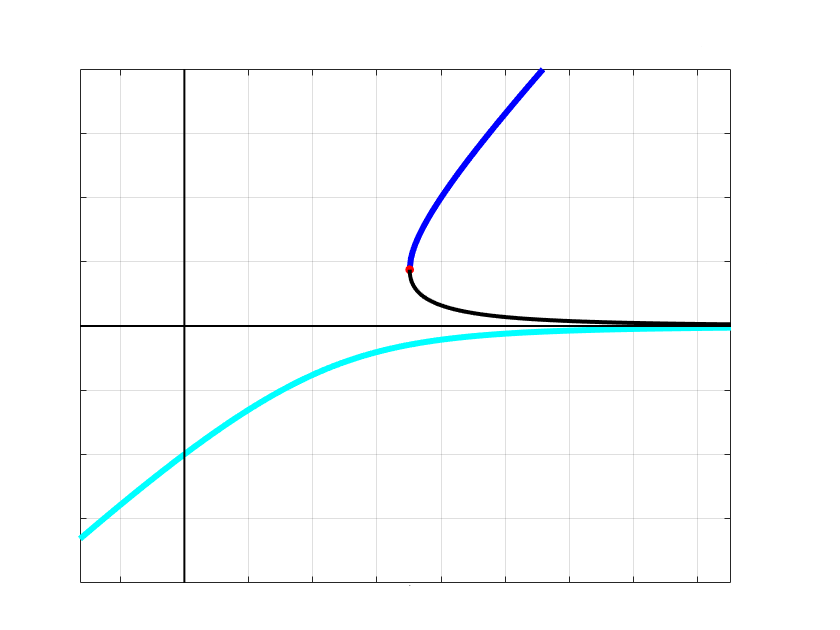}
\end{figure}

\begin{center}
  Plot of the set $\mathcal{C}$ for $\kappa = 1, \alpha = \frac{1}{3}$. The
  two black axes are $\{ z = 0 \}$ for the horiztonal one and $\{ y = 0 \}$
  for the vertical one.
\end{center}

We will show in section \ref{sec41} that $\mathcal{C}$ has this shape for any
$\kappa > 0, \alpha \in] 0, 1 [$. In particular, there exists a point $A$ (the
red dot) where a cusp happens on the curve, and its coordinates are
\[ A = (g (y_0), y_0) \]
with $y_0 \assign (\kappa \alpha)^{\frac{1}{1 + \alpha}}$. Remark that $g
(y_0) = \kappa^{\frac{1}{1 + \alpha}} \left( \alpha^{\frac{1}{1 + \alpha}} +
\alpha^{- \frac{\alpha}{1 + \alpha}} \right) > 0$. By the implicit function
theorem, we construct the functions
\[ y_+^{\ast} :] g (y_0), + \infty [\rightarrow] y_0, + \infty [(\tmop{in}
   \tmop{dark} \tmop{blue}), y^{\ast}_- : \mathbb{R} \rightarrow] - \infty, 0
   [ (\tmop{in} \tmop{cyan}) \]
as the respective inverses of $g :] y_0, + \infty [\rightarrow] g (y_0), +
\infty [$ and $g :] - \infty, 0 [\rightarrow \mathbb{R}$. We can now describe
the first order of the solution.

\begin{theorem}
  \label{thfo}For $\kappa > 0, \alpha \in] 0, 1 [$, there exists $z_c > g
  (y_0)$ such that the following holds. Consider $f$ the solution to the
  viscous Burgers equation $\partial_t f - \partial_x^2 f + f \partial_x f =
  0$ for an initial condition $f_0 \in C^1 (\mathbb{R}, \mathbb{R})$ with
  \[ f_0 (y) = \frac{\kappa (1 + o_{y \rightarrow \pm \infty} (1))}{| y
     |^{\alpha}} \]
  and
  \[ f'_0 (y) = \frac{- \alpha \kappa (1 + o_{y \rightarrow \pm \infty}
     (1))}{y | y |^{\alpha}} . \]
  Then, for any $z \in \mathbb{R} \backslash \{ z_c \}$, we have the
  convergence
  \[ t^{\frac{\alpha}{1 + \alpha}} f \left( z t^{\frac{1}{1 + \alpha}}, t
     \right) \rightarrow \mathfrak{p} (z) \]
  when $t \rightarrow + \infty$, where the profil $\mathfrak{p}$ is defined by
  \[ \mathfrak{p} (z) \assign \left\{\begin{array}{l}
       \kappa | y^{\ast}_+ (z) |^{- \alpha} \tmop{if} z > z_c\\
       \kappa | y_-^{\ast} (z) |^{- \alpha} \tmop{if} z < z_c .
     \end{array}\right. \]
  Furthermore, for any $\varepsilon > 0$, the convergence is uniform on
  $\mathbb{R} \backslash [z_c - \varepsilon, z_c + \varepsilon]$, and
  \[ \left\| t^{\frac{\alpha}{1 + \alpha}} f \left( . t^{\frac{1}{1 +
     \alpha}}, t \right) -\mathfrak{p} \right\|_{L^{\infty} (\mathbb{R}
     \backslash [z_c - \varepsilon, z_c + \varepsilon])} = O_{t \rightarrow +
     \infty} \left( t^{- \frac{1 - \alpha}{2 (1 + \alpha)}} \right) . \]
\end{theorem}

\begin{figure}[H]
    \centering
    \includegraphics[width=13cm]{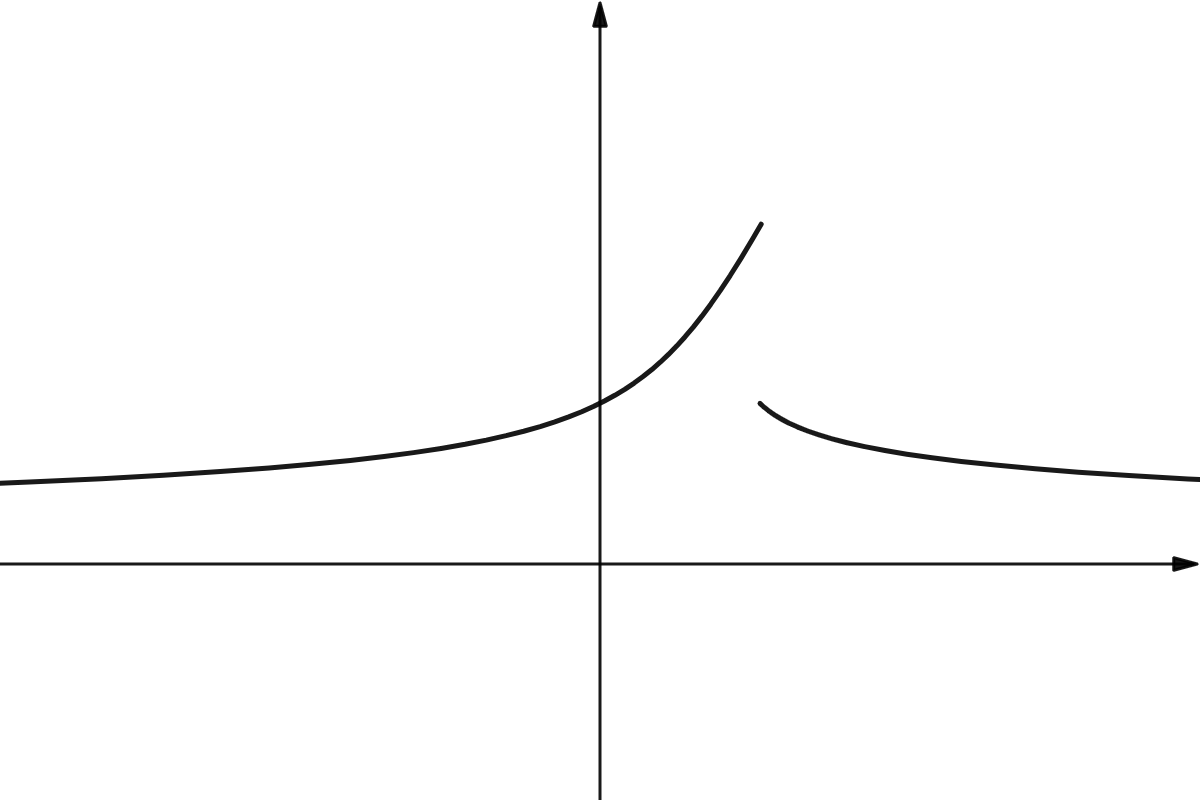}
\end{figure}

\begin{center}
  Graph of the function $z \rightarrow \mathfrak{p} (z)$. The point of
  discontinuity is $z_c$.
\end{center}

\

We can show that the profil $\mathfrak{p}$ is never continuous at $z = z_c$.
The computation of $z_c$ is rather intricate. We have it by an implicit
equation: it is the unique solution of $\mathcal{H} (y_+^{\ast} (z_c))
=\mathcal{H} (y_-^{\ast} (z_c))$ where
\[ \mathcal{H} (y) \assign - \frac{\kappa^2}{4 | y |^{2 \alpha}} -
   \frac{\kappa (1 - \alpha)}{2} | y |^{1 - \alpha} \]
in the set $] g (y_0), + \infty [$.

We can check that $y_+^{\ast} (z) - z \rightarrow 0$ when $z \rightarrow +
\infty$ and $y_-^{\ast} (z) - z \rightarrow 0$ when $z \rightarrow - \infty$.
This means that if $| z | \gg 1 + t$, we have $\mathfrak{p} (z) \simeq
\frac{\kappa}{| z |^{\alpha}} \simeq t^{\frac{\alpha}{1 + \alpha}} f_0 \left(
z t^{\frac{1}{1 + \alpha}} \right)$.

Remark that if $f (x, t)$ is solution of the viscous Burgers equation, then so
is $- f (- x, t)$. We can therefore consider negative initial data as well in
Theorems \ref{mainth} and \ref{thfo}.

This result does not say what is happening at $z_c$, but away from it, the
profile $\mathfrak{p}$ is a global attractor. However, from {\cite{GMP2}}, we
know that there are no generic profile near $z_c$. A consequence of Theorem
1.4 there is that in a vincinity of $z_c$, we can construct functions $f_0$
such that $t^{\frac{\alpha}{1 + \alpha}} f \left( . t^{\frac{1}{1 + \alpha}},
t \right)$ converges to different profiles near $z_c$.

\

As previously, it is interesting to compare this result with the first order
for the heat equation for the same initial data.

\begin{proposition}
  \label{heatfo}For $\kappa > 0, \alpha \in] 0, 1 [$, consider $f$ the
  solution of the heat equation $\partial_t f - \partial_x^2 f = 0$ for an
  initial condition
  \[ f_0 (y) = \frac{\kappa (1 + o_{| y | \rightarrow + \infty} (1))}{(1 + | y
     |)^{\alpha}} . \]
  Then, uniformly in $z \in \mathbb{R}$, we have the convergence
  \[ t^{\alpha / 2} f \left( \sqrt{t} z, t \right) \rightarrow
     \frac{\kappa}{\sqrt{4 \pi}} \int_{\mathbb{R}} \frac{1}{| y |^{\alpha}}
     e^{- (z - y)^2 / 4} d y \]
  when $t \rightarrow + \infty$.
\end{proposition}

This result has been shown in {\cite{GMP2}} (see Proposition 1.1 there).
Remark in particular that here, the limit profile is continuous. This also
means that the long time behaviour is highly nonlinear, as for instance has
been shown in {\cite{MR1266100}} in another context.

\subsection{Other exemples of discontinuous asymptotic profiles}

In this section, we compute the first order of solutions to the viscous
Burgers equation, but for other non integrable initial datas. Our first
exemple is the case where
\[ f_0 (y) = \frac{\mp \kappa (1 + o_{y \rightarrow \pm \infty} (1))}{| y
   |^{\alpha}}, \]
where we get a similar profile.

\begin{proposition}
  \label{3rdcase}For $\kappa > 0, \alpha \in] 0, 1 [$, there exists $z_d > 0$
  such that the following holds. Consider $f$ the solution to the viscous
  Burgers equation $\partial_t f - \partial_x^2 f + f \partial_x f = 0$ for an
  initial condition $f_0 \in C^1 (\mathbb{R}, \mathbb{R})$ with
  \[ f_0 (y) = \frac{- \kappa (1 + o_{y \rightarrow + \infty} (1))}{| y
     |^{\alpha}}, f_0 (y) = \frac{\kappa (1 + o_{y \rightarrow - \infty}
     (1))}{| y |^{\alpha}} \]
  and
  \[ f'_0 (y) = \frac{\alpha \kappa (1 + o_{y \rightarrow + \infty} (1))}{y |
     y |^{\alpha}}, f_0 (y) = \frac{- \alpha \kappa (1 + o_{y \rightarrow -
     \infty} (1))}{y | y |^{\alpha}} . \]
  Then, for any $z \in \mathbb{R} \backslash \{ z_d \}$, we have the
  convergence
  \[ t^{\frac{\alpha}{1 + \alpha}} f \left( z t^{\frac{1}{1 + \alpha}}, t
     \right) \rightarrow \mathfrak{p} (z) \]
  when $t \rightarrow + \infty$ to some profile $\mathfrak{p}$. Furthermore,
  for any $\varepsilon > 0$, the convergence is uniform on $\mathbb{R}
  \backslash [z_d - \varepsilon, z_d + \varepsilon]$, and
  \[ \left\| t^{\frac{\alpha}{1 + \alpha}} f \left( . t^{\frac{1}{1 +
     \alpha}}, t \right) -\mathfrak{p} \right\|_{L^{\infty} (\mathbb{R}
     \backslash [z_d - \varepsilon, z_d + \varepsilon])} = O_{t \rightarrow +
     \infty} \left( t^{- \frac{1 - \alpha}{2 (1 + \alpha)}} \right) . \]
\end{proposition}

The profile $\mathfrak{p}$ will have a similar definition as the one in
Theorem \ref{thfo}, and presents a discontinuity at some point $z_d$, but they
are not identical. We can check that $\mathfrak{p}$ is the unique entropic
solution in the sense of section 1.1.3 of {\cite{GMP2}} that behaves like
$\frac{\mp \kappa}{| y |^{\alpha}}$ at $\pm \infty$.

\

Interestingly, we cannot show a similar result for the last remaining case,
that is
\[ f_0 (y) = \frac{\pm \kappa (1 + o_{y \rightarrow \pm \infty} (1))}{| y
   |^{\alpha}}, \]
as our approach for the proof fails in that case. We believe that if the
convergences occurs, the profile will be of different type than the previous
ones (in particular, it will have at least two disconinuities). See subection
\ref{ss331} for more details about these two cases, as well as subsection
1.1.3 of {\cite{GMP2}}.

\

We now consider other type of decay, namely the case
\[ f_0 (y) = \frac{\kappa (1 + o_{y \rightarrow \pm \infty} (1))}{| y
   |^{\alpha} \ln^{\beta} (| y |)} \]
for $\kappa, \beta > 0, \alpha \in] 0, 1 [$. It turns out that the functions
$y_{\pm}^{\ast}$ defined for Theorem \ref{thfo} will also appear here.

\begin{proposition}
  \label{lng}For $\kappa > 0, \alpha \in] 0, 1 [, \beta > 0$, there exists
  $z_c > g (y_0)$ such that the following holds. Consider $f$ the solution to
  the viscous Burgers equation $\partial_t f - \partial_x^2 f + f \partial_x f
  = 0$ for an initial condition $f_0 \in C^1 (\mathbb{R}, \mathbb{R})$ with
  \[ f_0 (y) = \frac{\kappa (1 + o_{y \rightarrow \pm \infty} (1))}{| y
     |^{\alpha} \ln^{\beta} (| y |)} \]
  and
  \[ f'_0 (y) = \frac{- \alpha \kappa (1 + o_{y \rightarrow \pm \infty}
     (1))}{y | y |^{\alpha} \ln^{\beta} (| y |)} . \]
  Then, for any $z \in \mathbb{R} \backslash \{ z_c \}$, we have, with $\mu
  (t)$ the solution for $t$ large enough to
  \[ \mu (t)^{1 + \alpha} \ln^{\beta} (\mu (t)) = t, \]
  the convergence
  \[ \frac{t}{\mu (t)} f (z \mu (t), t) \rightarrow \mathfrak{p} (z) \]
  when $t \rightarrow + \infty$, where the profil $\mathfrak{p}$ is the same
  as Theorem \ref{thfo}. Furthermore, for any $\varepsilon > 0$, the
  convergence is uniform on $\left[ - \frac{1}{\varepsilon},
  \frac{1}{\varepsilon} \right] \backslash [z_c - \varepsilon, z_c +
  \varepsilon]$.
\end{proposition}

Here, the limite profile $\mathfrak{p}$ is the same as in the case $\beta = 0$
(that is Theorem \ref{thfo}), but the scalings are different. A similar result
can be proven for a larger class of asymptotics on the initial condition, but
writing a general result is still difficult at this point.

Remark that at fixed $\alpha$ but for different values of $\beta$, the initial
datas in Proposition \ref{lng} are in the same $L^p (\mathbb{R})$ space, and
not in the same $L^q (\mathbb{R})$, but the scaling $\mu (t)$ is different.
This means that we can not expect the asymptotique profile to only depends in
which $L^p (\mathbb{R})$ spaces the initial data is or is not to deduce the
asymptotic profile, we need more informations than that.

\

More generically, our approach should work for initial datas $f_0$ such that
there exists $\mu (t)$ and $\mathfrak{f}$ a $C^1$ function that never cancels
such that
\[ \frac{t}{\mu (t)} f_0 (\mu (t) y) \rightarrow \mathfrak{f} (y), t f'_0 (\mu
   (t) z) \rightarrow \mathfrak{f}' (y) \]
for $y \neq 0$ when $t \rightarrow + \infty$, and with $y \rightarrow y
+\mathfrak{f} (y) : \mathbb{R}^{\ast} \rightarrow \mathbb{R}$ surjective. The
asymptotic profile will then depend only on $\mathfrak{f}$. Remark that for
either $f_0 (x) \sim \frac{\kappa}{| x |^{\alpha}}$ and $\mu (t) =
t^{\frac{1}{1 + \alpha}}$ or $f_0 (x) \sim \frac{\kappa}{| y |^{\alpha}
\ln^{\beta} (| y |)}$ and $\mu (t)$ solution of $\mu (t)^{1 + \alpha}
\ln^{\beta} (\mu (t)) = t$, we have
\[ \frac{t}{\mu (t)} f_0 (\mu (t) y) \rightarrow \frac{\kappa}{| y |^{\alpha}}
\]
for $y \neq 0$, and this is why both profiles are identical, only the scaling
has changed.

\

Finally, as a last exemple, we look at the non symmetric case, when the decay
rate is not the same at $\pm \infty$. We require there more precisions on the
behavior of $f_0$. This is because we are in a case where $y +\mathfrak{f}
(y)$ will not be surjective.

\begin{proposition}
  \label{IDwhynot}For $\kappa > 0, \alpha \in] 0, 1 [$ there exists $z_e > 0$
  such that the following holds. Given $1 > \beta > \alpha$, consider $f$ the
  solution to the viscous Burgers equation $\partial_t f - \partial_x^2 f + f
  \partial_x f = 0$ for an initial condition $f_0 \in C^1 (\mathbb{R},
  \mathbb{R})$ with
  \[ f_0 (y) = \frac{\kappa (1 + o_{y \rightarrow + \infty} (1))}{| y
     |^{\alpha}}, f_0 (y) = \frac{\kappa (1 + o_{y \rightarrow - \infty}
     (1))}{| y |^{\beta}} \]
  and
  \[ f'_0 (y) = \frac{- \alpha \kappa (1 + o_{y \rightarrow + \infty} (1))}{y
     | y |^{\alpha}}, f_0 (y) = \frac{- \beta \kappa (1 + o_{y \rightarrow -
     \infty} (1))}{y | y |^{\beta}} . \]
  Suppose also that $f_0 > 0$ and $f_0' > 0$ on $] - \infty, 0 [$ and $f_0' <
  0$ on $] 0, + \infty [$. Then, for any $z \in \mathbb{R} \backslash \{ 0,
  z_e \}$, we have the convergence
  \[ t^{\frac{\alpha}{1 + \alpha}} f \left( z t^{\frac{1}{1 + \alpha}}, t
     \right) \rightarrow \mathfrak{p} (z) \]
  when $t \rightarrow + \infty$, where the profil $\mathfrak{p}$ is defined by
  \[ \mathfrak{p} (z) \assign \left\{\begin{array}{l}
       \kappa | y^{\ast}_+ (z) |^{- \alpha} \tmop{if} z > z_e\\
       z \tmop{if} 0 < z < z_e\\
       0 \tmop{if} z < 0.
     \end{array}\right. \]
\end{proposition}

Remark that here also, the profile has a discontinuity at $z_e$, and its
derivative also has a discontinuity at $0$. Here, the branch $y^{\ast}_+$ is
the same as in Theorem \ref{thfo}, but for the same value of $\alpha$ and
$\kappa$, we do not have $z_e = z_c$. The fact that the identity appears on
some interval is consistent with the entropic arguments of section 1.1.3 of
{\cite{GMP2}}.

\

In all these exemples, we have taken the same value of $\kappa$ at $+ \infty$
and $- \infty$, but in fact this is not necessary for our proofs. It is taken
to simplify some computations and notations, but similar results should hold
in the more general setting. See subsection \ref{54ss13} for the main steps of
the proof of these results.

\subsection{Plan of the proofs\label{54ss13}}

The proofs of Theorem \ref{mainth} and Proposition \ref{Prop14} rely on the
Hopf-Cole transformation ({\cite{MR42889}}, {\cite{MR47234}}), and they are
done in section \ref{sec2}. By this transformation we show that the solution
of the viscous Burgers equation is
\[ f (x, t) = \frac{\int_{\mathbb{R}} \frac{x - y}{t} \exp \left( - \frac{(x -
   y)^2}{4 t} - \frac{1}{2} \int_0^y f_0 (z) d z \right) d
   y}{\int_{\mathbb{R}} \exp \left( - \frac{(x - y)^2}{4 t} - \frac{1}{2}
   \int_0^y f_0 (z) d z \right) d y} . \]
We check that this quantity makes sense even with the slow decay of $f_0$.
After an integration by parts, we write it
\begin{equation}
  f (x, t) = \frac{\int_{\mathbb{R}} f_0 (y) e^{H (y)} d y}{\int_{\mathbb{R}}
  e^{H (y)} d y} \label{eq11}
\end{equation}
where $H (y) \assign - \frac{(x - y)^2}{4 t} - \frac{1}{2} \int_0^y f_0 (z) d
z$. The core of our proof is to show that, although it is difficult to
estimate $\int_{\mathbb{R}} e^{H (y)} d y$, most of its mass (up to an almost
exponentially small error in time) comes from the integral on $\mathbb{R}
\backslash \left[ - \varepsilon t^{\frac{1}{1 + \alpha}}, \varepsilon
t^{\frac{1}{1 + \alpha}} \right]$ for some small $\varepsilon > 0$. Therefore,
since $| f_0 (y) |$ decays like $| y |^{- \alpha}$ and does not cancel, we
have
\[ \int_{\mathbb{R}} | f_0 (y) | e^{H (y)} d y \lesssim \left(
   \sup_{\mathbb{R} \backslash \left[ - \varepsilon t^{\frac{1}{1 + \alpha}},
   \varepsilon t^{\frac{1}{1 + \alpha}} \right]} | f_0 | \right)
   \int_{\mathbb{R}} e^{H (y)} d y \leqslant \frac{K}{t^{\frac{\alpha}{1 +
   \alpha}}} \int_{\mathbb{R}} e^{H (y)} d y, \]
leading to the estimate on $f$. To show Proposition \ref{Prop14}, we
differentiate equation (\ref{eq11}) and we check that $\partial_t^n
\partial_x^k f$ can be estimated by terms of the form $\frac{\int_{\mathbb{R}}
g_0 (y) e^{H (y)} d y}{\int_{\mathbb{R}} e^{H (y)} d y}$ where $g_0$ is a sum
of derivatives and powers of $f_0$. The proof follows from similar arguments
as the one of Theorem \ref{mainth}.

\

The proof of Theorem \ref{thfo} follows similar ideas, but with the
additional precision on $f_0$ at infinity, we can compute exactly where $H
\left( .t^{\frac{1}{1 + \alpha}} \right)$ reaches its maximum. It turns out
that at most two particular functions of $z$, that we will call $y_- (z, t)$
and $y_+ (z, t)$, can reach it. When $t \rightarrow + \infty$ they both
converge nicely to $y_+^{\ast} (z)$ and $y_-^{\ast} (z)$ respectively, at
least where they have a chance to reach the maximum. For $z$ close to $-
\infty$ it will be reached only by $y_-$, and close to $+ \infty$ only by
$y_+$. We show that there exists only one value of $z$ (which is $z_c$) such
that the maximum is reached by both. Then, if $z \neq z_c$, most of the mass
of $\int_{\mathbb{R}} e^{H (y)} d y$ is coming from a small neighborhood of
the point where $H \left( .t^{\frac{1}{1 + \alpha}} \right)$ reaches its
maximum, and thus
\[ \int_{\mathbb{R}} f_0 (y) e^{H (y)} d y \simeq f_0 \left( y_{\max}
   t^{\frac{1}{1 + \alpha}} \right) \int_{\mathbb{R}} e^{H (y)} d y \simeq -
   \kappa t^{- \frac{\alpha}{1 + \alpha}} | y_{\max} |^{- \alpha}
   \int_{\mathbb{R}} e^{H (y)} d y \]
where $y_{\max}$ is the value reaching the maximum. This result is proven in
section \ref{sec4}. Proposition \ref{3rdcase} to \ref{IDwhynot} are proven
using similar arguments, see section \ref{ss33}.

\

Propositions \ref{prop12} and \ref{prop15} follow from computations on the
explicit solution of the heat equation $\partial_t f - \partial_x^2 f = 0$
(see {\cite{MR2597943}}):
\[ f (x, t) = \frac{1}{\sqrt{4 \pi t}} \int_{\mathbb{R}} f_0 (y) e^{- \frac{(x
   - y)^2}{4 t}} d y. \]
Although their proofs are not difficult, for the sake of completness, they are
done in annex \ref{sec3}.

\subsection{Open problems}

Our approach does not work if the initial data satisfies
\[ f_0 (y) = \frac{\pm \kappa (1 + o_{y \rightarrow \pm \infty} (1))}{| y
   |^{\alpha}} . \]
We can show that then, $t^{\frac{\alpha}{1 + \alpha}} f \left( z t^{\frac{1}{1
+ \alpha}}, t \right)$ converges to some profile $\mathfrak{p}$ outside of a
compact interval, but we have no informations on what is happening on it. It
is likely that the solution is unstable and depends also on $f_0$ and not
simply on its behavior at $\pm \infty$.

It would also be interesting to understand how much of the results is still
true if we consider initial datas satisfying
\[ \frac{\kappa_2}{(1 + | y |)^{\alpha}} \geqslant f_0 (y) \geqslant
   \frac{\kappa_1}{(1 + | y |)^{\beta}} \]
if $1 > \alpha > \beta > 0$, or for other type of behavior at infinity of
$f_0$.

\

Concerning similar problem on other equations, since the proofs here use the
Hopf-Cole transform, it is unlikely that they will hold in a more general
setting. In {\cite{GMP2}} we show the existence of a local attractor (instead
of a global one) for the equation $\partial_t u - \partial_x^2 u + \partial_x
\left( \frac{u^2}{2} + J (u) \right) = 0$ if $| J (u) | \leqslant C | u |^3$,
see Proposition 1.5 there. A natural generalisation would be to look at
similar equations in dimension 2 or higher.

\begin{acknowledgments*}
  The authors are supported by Tamkeen under the NYU Abu Dhabi Research
  Institute grant CG002. The authors have no competing interests to declare that are relevant to the content of this article.
\end{acknowledgments*}

\section{Estimates on the viscous Burgers equation}\label{sec2}

This section is devoted to the proofs of Theorem \ref{mainth} and Proposition
\ref{Prop14}. We recall that the quantity
\begin{equation}
  f (x, t) = \frac{\int_{\mathbb{R}} \frac{x - y}{t} \exp \left( - \frac{(x -
  y)^2}{4 t} - \frac{1}{2} \int_0^y f_0 (z) d z \right) d y}{\int_{\mathbb{R}}
  \exp \left( - \frac{(x - y)^2}{4 t} - \frac{1}{2} \int_0^y f_0 (z) d z
  \right) d y}
\end{equation}
is well defined for all $t \geqslant 0, x \in \mathbb{R}$ and is the solution
to the viscous burgers equation $\partial_t f - \partial_x^2 f + f \partial_x
f = 0$ with initial data $f_0$. We compute, by integration by parts, that
\begin{eqnarray*}
  &  & \int_{\mathbb{R}} \frac{x - y}{t} \exp \left( - \frac{(x - y)^2}{4 t}
  - \frac{1}{2} \int_0^y f_0 (z) d z \right) d y\\
  & = & \int_{\mathbb{R}} 2 \partial_y \left( \exp \left( - \frac{(x -
  y)^2}{4 t} \right) \right) \exp \left( - \frac{1}{2} \int_0^y f_0 (z) d z
  \right) d y\\
  & = & \int_{\mathbb{R}} f_0 (y) \exp \left( - \frac{(x - y)^2}{4 t} -
  \frac{1}{2} \int_0^y f_0 (z) d z \right) d y.
\end{eqnarray*}
This integration by part is justified for any $t > 0$ and $x \in \mathbb{R}$
by the fact that
\[ \left| \int_0^y f_0 (z) d z \right| \leqslant K | y |^{1 - \alpha}
   \leqslant \frac{(x - y)^2}{8 t} \]
for $| y |$ large enough. We define the quantity
\begin{equation}
  H (y) = - \frac{(x - y)^2}{4 t} - \frac{1}{2} \int_0^y f_0 (z) d z,
\end{equation}
so that
\begin{equation}
  f (x, t) = \frac{\int_{\mathbb{R}} f_0 (y) e^{H (y)} d y}{\int_{\mathbb{R}}
  e^{H (y)} d y} .
\end{equation}
In particular we have the inequality
\[ | f (x, t) | \leqslant \| f_0 \|_{L^{\infty}} . \]
We study here the case $\frac{- \kappa_2}{(1 + | y |)^{\alpha}} \leqslant f_0
(y) \leqslant \frac{- \kappa_1}{(1 + | y |)^{\alpha}}$. If $f (x, t)$ solves
the viscous Burgers equation, then so does $- f (- x, t)$, so this is
equivalent as considering $\frac{\kappa_2}{(1 + | y |)^{\alpha}} \geqslant f_0
(y) \geqslant \frac{\kappa_1}{(1 + | y |)^{\alpha}}$. We do so because some
estimates will now be in a more usual direction. We consider the more general
quantity
\[ F (x, t) \assign \frac{\int_{\mathbb{R}} g (y) e^{H (y)} d
   y}{\int_{\mathbb{R}} e^{H (y)} d y}, \]
where we still take $H (y) = - \frac{(x - y)^2}{4 t} - \frac{1}{2} \int_0^y
f_0 (z) d z$. For now, we simply suppose that $g \in C^0 (\mathbb{R})$ with
$\| g \|_{C^0 (\mathbb{R})} \leqslant \kappa$ where $\kappa$ only depends on
$\kappa_1, \kappa_2$ and $\alpha$, so that in particular the quantity
$\int_{\mathbb{R}} g (y) e^{H (y)} d y$ is well defined. Let us first show the
following result.

\

\begin{proposition}
  \label{prop211}There exists $T^{\ast}, K, \nu > 0, \mu_1 < 0, \mu_2 > 0$
  depending on $\kappa_1, \kappa_2, \alpha$ such that, if we define the set $D
  \assign \left] - \infty, \mu_1 t^{\frac{1}{1 + \alpha}} \right] \cup \left[
  \mu_2 t^{\frac{1}{1 + \alpha}}, + \infty \right[$, then for $x \in
  \mathbb{R}, t \geqslant T^{\ast},$
  \[ | F (x, t) | \leqslant \sup_D | g | + K e^{- \nu t^{\frac{1 - \alpha}{1 +
     \alpha}}} . \]
  \[ \  \]
\end{proposition}

\subsection{Proof of Proposition \ref{prop211}}

We first estimate
\[ | F (x, t) | \leqslant \sup_D | g | \frac{\int_D e^{H (y)} d
   y}{\int_{\mathbb{R}} e^{H (y)} d y} + \| g \|_{L^{\infty}}
   \frac{\int_{\mathbb{R} \backslash D} e^{H (y)} d y}{\int_{\mathbb{R}} e^{H
   (y)} d y} \leqslant \sup_D | g | + \kappa, \]
therefore the result hold on $[0, T^{\ast}]$ for any constant $T^{\ast} > 0$
depending on $\kappa_1, \kappa_2, \alpha$.

We now suppose that $t \geqslant T^{\ast}$.

\

We compute
\[ H' (y) = \frac{1}{2 t} (x - y) - \frac{1}{2} f_0 (y) = \frac{1}{2} \left(
   \frac{x - y}{t} - f_0 (y) \right) . \]
Since $H (y) \rightarrow - \infty$ when $y \rightarrow \pm \infty$, the
function $H$ must reach its maximum at a point where $H' (y) = 0$. We consider
any such $y^{\ast} \in \mathbb{R}$ solution of $H' (y^{\ast}) = 0.$ Let us
give some estimates on $y^{\ast}$. We have
\[ \frac{1}{2 t} (x - y^{\ast}) + \frac{\kappa_1}{2 (1 + | y^{\ast}
   |)^{\alpha}} \leqslant \frac{1}{2 t} (x - y^{\ast}) - \frac{1}{2} f_0
   (y^{\ast}) = 0 \leqslant \frac{1}{2 t} (x - y^{\ast}) + \frac{\kappa_2}{2
   (1 + | y^{\ast} |)^{\alpha}}, \]
therefore
\begin{equation}
  - \kappa_2 t \leqslant (1 + | y^{\ast} |)^{\alpha} (x - y^{\ast}) \leqslant
  - \kappa_1 t. \label{agh}
\end{equation}
We will use different arguments to show the estimate on $F$ depending on the
value of $\frac{x}{t^{\frac{1}{1 + \alpha}}}$. We decompose the problem into
four different cases.

\

In the rest of the proof, we will use $\varepsilon > 0$ a small constant
depending on $\kappa_1, \kappa_2, \alpha$, and then we will take $T^{\ast}$
large, depending on $\kappa_1, \kappa_2, \alpha$ and $\varepsilon$. In this
section, a generic constant $K$ can depend on $\kappa_1, \kappa_2$ and
$\alpha$, but not on $\varepsilon$ or $T^{\ast}$, except if it is explicitely
stated.

\

\subsubsection{Estimate on $F$ in the case $\frac{| x |}{t^{\frac{1}{1 +
\alpha}}} \leqslant \varepsilon$}

We recall that we take $t \geqslant T^{\ast} > 0$ with $T^{\ast}$ large.
Define $z^{\ast} = \frac{y^{\ast}}{t^{\frac{1}{1 + \alpha}}}$. By (\ref{agh}),
it satisfies
\begin{equation}
  - \kappa_2 \leqslant \left( \frac{1}{t^{\frac{1}{1 + \alpha}}} + | z^{\ast}
  | \right)^{\alpha} \left( \frac{x}{t^{\frac{1}{1 + \alpha}}} - z^{\ast}
  \right) \leqslant - \kappa_1 . \label{654987}
\end{equation}
First, let us show that for $t$ large enough, we have $z^{\ast} >
\varepsilon$. Indeed, if it's not true, then $z^{\ast} \leqslant \varepsilon$
and by (\ref{654987}) since
\[ \frac{x}{t^{\frac{1}{1 + \alpha}}} - z^{\ast} \leqslant 0, \]
we have $z^{\ast} \geqslant \frac{x}{t^{\frac{1}{1 + \alpha}}}$, therefore $|
z^{\ast} | \leqslant \max \left( \varepsilon, \frac{| x |}{t^{\frac{1}{1 +
\alpha}}} \right) \leqslant \varepsilon$. But then,

\[ \left| \left( \frac{1}{t^{\frac{1}{1 + \alpha}}} + | z^{\ast} |
   \right)^{\alpha} \left( \frac{x}{t^{\frac{1}{1 + \alpha}}} - z^{\ast}
   \right) \right| \leqslant \left( \frac{1}{t^{\frac{1}{1 + \alpha}}} +
   \varepsilon \right)^{\alpha} 2 \varepsilon < \kappa_1 \]
for $\varepsilon$ small enough and $T^{\ast}$ large enough (depending on
$\varepsilon, \kappa_1, \kappa_2$), leading to a contradiction.

\

Secondly, there exists $\kappa > 1$ large enough (depending only on
$\kappa_2, \varepsilon$) such that
\[ z^{\ast} \leqslant \kappa . \]
Indeed, otherwise $z^{\ast} \geqslant \kappa$ and $\frac{x}{t^{\frac{1}{1 +
\alpha}}} - z^{\ast} \leqslant \varepsilon - \kappa < 0$, leading to
\[ - \kappa_2 \leqslant \left( \frac{1}{t^{\frac{1}{1 + \alpha}}} + | z^{\ast}
   | \right)^{\alpha} \left( \frac{x}{t^{\frac{1}{1 + \alpha}}} - z^{\ast}
   \right) \leqslant (2 \kappa)^{\alpha} (- \kappa + \varepsilon) < - \kappa_2
\]
if $\kappa$ and $T^{\ast}$ are large enough.

We summarize. There exists $\Lambda_2 > \Lambda_1 > 0$ two constants depending
on $\kappa_1, \kappa_2, \varepsilon$ such that if $\frac{| x |}{t^{\frac{1}{1
+ \alpha}}} \leqslant \varepsilon$, then any $y^{\ast} \in \mathbb{R}$
solution of $H' (y^{\ast}) = 0$ is in the set $\left[ \Lambda_1 t^{\frac{1}{1
+ \alpha}}, \Lambda_2 t^{\frac{1}{1 + \alpha}} \right]$. Since $H (y)
\rightarrow - \infty$ when $y \rightarrow \pm \infty$, its maximum is reached
in this set. For such a $y^{\ast} \in \mathbb{R}$, since $f_0 (y) \leqslant
\frac{- \kappa_1}{(1 + | y |)^{\alpha}}$ we have
\[ - \frac{1}{2} \int_0^{y^{\ast}} f_0 (z) d z \geqslant \frac{\kappa_1}{4 (1
   - \alpha)} (1 + | y^{\ast} |)^{1 - \alpha} \geqslant \frac{\kappa_1
   \Lambda_1}{4 (1 - \alpha)} t^{\frac{1 - \alpha}{1 + \alpha}} \]
and
\[ (x - y^{\ast})^2 \leqslant (\varepsilon - \Lambda_2)^2 t^{\frac{2}{1 +
   \alpha}}, \]
leading to
\[ - \frac{(x - y^{\ast})^2}{4 t} \geqslant - \frac{(\varepsilon -
   \Lambda_2)^2}{4} t^{\frac{1 - \alpha}{1 + \alpha}} . \]
We deduce that
\[ H (y^{\ast}) \geqslant \left( \frac{\kappa_1 \Lambda_1}{4 (1 - \alpha)} -
   \frac{(\varepsilon - \Lambda_2)^2}{4} \right) t^{\frac{1 - \alpha}{1 +
   \alpha}} . \]
We define
\[ C_0 \assign t^{- \frac{1 - \alpha}{1 + \alpha}} \max_{\mathbb{R}} H \]
which therefore satisfies $C_0 \geqslant \frac{\kappa_1 \Lambda_2}{4 (1 -
\alpha)} - \frac{(\varepsilon - \Lambda_1)^2}{4}$. Remark that $C_0$ depends
on time and it is possible that $C_0 = 0$ or $C_0 < 0$. We check similarly
that $C_0 \leqslant K (\varepsilon)$, thus $C_0$ is uniformly bounded in time.

\

Now, we estimate for $\gamma > 0$ that
\[ H' \left( \gamma t^{\frac{1}{1 + \alpha}} \right) \geqslant
   \frac{t^{\frac{- \alpha}{1 + \alpha}}}{2} \left( \frac{x}{t^{\frac{1}{1 +
   \alpha}}} - \gamma + \frac{\kappa_1}{\left( t^{- \frac{1}{1 + \alpha}} +
   \gamma \right)^{\alpha}} \right) . \]
We deduce that if $T^{\ast}$ is large enough, there exists $\lambda_1 > 0$
with $\lambda_1 < \Lambda_1$ such that $\frac{x}{t^{\frac{1}{1 + \alpha}}} -
\gamma + \frac{\kappa_1}{\left( t^{- \frac{1}{1 + \alpha}} + \gamma
\right)^{\alpha}} \geqslant 2$ for $\gamma \in \left[ \frac{\lambda_1}{2},
\lambda_1 \right]$ and thus $H' (y) \geqslant t^{\frac{- \alpha}{1 + \alpha}}$
for $y \in \left[ \frac{\lambda_1}{2} t^{\frac{1}{1 + \alpha}}, \lambda_1
t^{\frac{1}{1 + \alpha}} \right]$. Similarly, we check that for $\gamma > 0$,
\[ H' \left( \gamma t^{\frac{1}{1 + \alpha}} \right) \leqslant
   \frac{t^{\frac{- \alpha}{1 + \alpha}}}{2} \left( \frac{x}{t^{\frac{1}{1 +
   \alpha}}} - \gamma + \frac{\kappa_2}{\left( t^{- \frac{1}{1 + \alpha}} +
   \gamma \right)^{\alpha}} \right), \]
and that there exists $\lambda_2 > \Lambda_2$ such that
$\frac{x}{t^{\frac{1}{1 + \alpha}}} - \gamma + \frac{\kappa_2}{\left( t^{-
\frac{1}{1 + \alpha}} + \gamma \right)^{\alpha}} \leqslant - 2$ if $\gamma >
\lambda_2$. This implies that $H' (y) \leqslant - t^{\frac{- \alpha}{1 +
\alpha}}$ for $y \geqslant \lambda_2 t^{\frac{1}{1 + \alpha}}$.

Since $H'$ only has zeros in $\left[ \Lambda_1 t^{\frac{1}{1 + \alpha}},
\Lambda_2 t^{\frac{1}{1 + \alpha}} \right]$ and $H (y) \rightarrow - \infty$
when $y \rightarrow - \infty$, we have that $H' \geqslant 0$ on $\left] -
\infty, \Lambda_1 t^{\frac{1}{1 + \alpha}} \right[$. Therefore, for $y \in
\left] - \infty, \frac{\lambda_1}{2} t^{\frac{1}{1 + \alpha}} \right[$, we
have
\begin{eqnarray*}
  H (y) & = & H \left( \lambda_1 t^{\frac{1}{1 + \alpha}} \right) -
  \int_y^{\lambda_1 t^{\frac{1}{1 + \alpha}}} H' (z) d z\\
  & \leqslant & \max H - \int_{\frac{\lambda_1}{2} t^{\frac{1}{1 +
  \alpha}}}^{\lambda_1 t^{\frac{1}{1 + \alpha}}} H' (z) d z\\
  & \leqslant & C_0 t^{\frac{1 - \alpha}{1 + \alpha}} - \frac{\lambda_1}{2}
  t^{\frac{1}{1 + \alpha} + \frac{- \alpha}{1 + \alpha}}\\
  & \leqslant & \left( C_0 - \frac{\lambda_1}{2} \right) t^{\frac{1 -
  \alpha}{1 + \alpha}} .
\end{eqnarray*}
Furthermore, for $y \in \left] 2 \lambda_2 t^{\frac{1}{1 + \alpha}}, + \infty
\right[$ we have
\begin{eqnarray*}
  H (y) & \leqslant & \max H + \int_{\lambda_2 t^{\frac{1}{1 + \alpha}}}^{2
  \lambda_2 t^{\frac{1}{1 + \alpha}}} H' (z) d z\\
  & \leqslant & (C_0 - \lambda_2) t^{\frac{1 - \alpha}{1 + \alpha}} .
\end{eqnarray*}
We define the domain
\[ D_1 \assign \left[  \frac{\lambda_1}{2} t^{\frac{1}{1 + \alpha}}, 2
   \lambda_2 t^{\frac{1}{1 + \alpha}} \right] . \]
We have shown that there exists $1 > \nu > 0$ a small constant depending on
$\kappa_1, \kappa_2, \varepsilon, T^{\ast}$ such that, for $y \in \mathbb{R}
\backslash D_1$, $H (y) \leqslant (C_0 - \nu) t^{\frac{1 - \alpha}{1 +
\alpha}} .$ Now we estimate
\[ \left| \int_{\mathbb{R}} g (y) e^{H (y)} d y \right| \leqslant \max_{D_1}
   (| g |) \int_{D_1} e^{H (y)} d y + \| g \|_{L^{\infty} (\mathbb{R})}
   \int_{\mathbb{R} \backslash D_1} e^{H (y)} d y \]
and thus
\[ | F (x, t) | \leqslant \max_{D_1} (| g |) + K \frac{\int_{\mathbb{R}
   \backslash D_1} e^{H (y)} d y}{\int_{\mathbb{R}} e^{H (y)} d y} . \]
Since $H' (y) = \frac{1}{2} \left( \frac{x - y}{t} - f_0 (y) \right)$, for $|
x | \leqslant \varepsilon t^{\frac{1}{1 + \alpha}}$ and $y \in \left[
\Lambda_1 t^{\frac{1}{1 + \alpha}}, \Lambda_2 t^{\frac{1}{1 + \alpha}}
\right]$, if $T^{\ast}$ is large enough we have $| H' (y) | \leqslant K$.
Therefore, there exists $\rho > 0$ depending only on $\kappa_1, \kappa_2$ such
that, for a point $y^{\ast}$ with $H (y^{\ast}) = \max H$ and $y \in [y^{\ast}
- \rho, y^{\ast} + \rho]$, we have $H (y) \geqslant (\max H) - 1$. We deduce
that
\[ \int_{\mathbb{R}} e^{H (y)} d y \geqslant \int_{[y^{\ast} - \rho, y^{\ast}
   + \rho]} e^{H (y)} d y \geqslant K (\kappa_1, \kappa_2) e^{\max H}
   \geqslant K (\kappa_1, \kappa_2) e^{C_0 t^{\frac{1 - \alpha}{1 + \alpha}}}
   . \]
We recall that
\[ H (y) = - \frac{(x - y)^2}{4 t} - \frac{1}{2} \int_0^y f_0 (z) d z, \]
and since
\[ \left| \int_0^y f_0 (z) d z \right| \leqslant K (1 + | y |)^{1 - \alpha},
\]
there exists $\Omega_1 < \frac{\lambda_1}{2}, \Omega_2 > 2 \lambda_2$ (with
eventually $\Omega_1 < 0$) depending on $\kappa_1, \kappa_2, \varepsilon,
T^{\ast}$, such that, outside of $\left[ \Omega_1 t^{\frac{1}{1 + \alpha}},
\Omega_2 t^{\frac{1}{1 + \alpha}} \right]$, we have
\begin{equation}
  H (y) \leqslant - \frac{(x - y)^2}{8 t} + (C_0 - \nu) t^{\frac{1 - \alpha}{1
  + \alpha}} \label{123}
\end{equation}
(we recall that here $\frac{| x |}{t^{\frac{1}{1 + \alpha}}} \leqslant
\varepsilon$ and $C_0$ is uniformly bounded in time). We deduce that
\begin{eqnarray*}
  \int_{\mathbb{R} \backslash \left[ \Omega_1 t^{\frac{1}{1 + \alpha}},
  \Omega_2 t^{\frac{1}{1 + \alpha}} \right]} e^{H (y)} d y & \leqslant &
  e^{(C_0 - \nu) t^{\frac{1 - \alpha}{1 + \alpha}}} \int_{\mathbb{R}
  \backslash \left[ \Omega_1 t^{\frac{1}{1 + \alpha}}, \Omega_2 t^{\frac{1}{1
  + \alpha}} \right]} e^{- \frac{(x - y)^2}{8 t}} d y\\
  & \leqslant & e^{(C_0 - \nu) t^{\frac{1 - \alpha}{1 + \alpha}}}
  \int_{\mathbb{R}} e^{- \frac{y^2}{8 t}} d y\\
  & \leqslant & K e^{(C_0 - \nu) t^{\frac{1 - \alpha}{1 + \alpha}}} \sqrt{t},
\end{eqnarray*}
and we compute
\begin{eqnarray*}
  &  & \int_{\left[ \Omega_1 t^{\frac{1}{1 + \alpha}}, \Omega_2 t^{\frac{1}{1
  + \alpha}} \right] \backslash D_1} e^{H (y)} d y\\
  & \leqslant & e^{(C_0 - \nu) t^{\frac{1 - \alpha}{1 + \alpha}}} \tmop{Vol}
  \left( \left[ \Omega_1 t^{\frac{1}{1 + \alpha}}, \Omega_2 t^{\frac{1}{1 +
  \alpha}} \right] \backslash D_1 \right)\\
  & \leqslant & K (\varepsilon, T^{\ast}) t^{\frac{1}{1 + \alpha}} e^{(C_0 -
  \nu) t^{\frac{1 - \alpha}{1 + \alpha}}} .
\end{eqnarray*}
We deduce that
\begin{eqnarray*}
  &  & \frac{\int_{\mathbb{R} \backslash D_1} e^{H (y)} d
  y}{\int_{\mathbb{R}} e^{H (y)} d y}\\
  & \leqslant & \frac{K (\varepsilon, T^{\ast}) \left(  \sqrt{t} e^{(C_0 -
  \nu) t^{\frac{1 - \alpha}{1 + \alpha}}} + t^{\frac{1}{1 + \alpha}} e^{(C_0 -
  \nu) t^{\frac{1 - \alpha}{1 + \alpha}}} \right)}{e^{C_0 t^{\frac{1}{1 +
  \alpha}}}}\\
  & \leqslant & K (\varepsilon, T^{\ast}) e^{- \nu t^{\frac{1 - \alpha}{1 +
  \alpha}}} \left( \sqrt{t} + t^{\frac{1}{1 + \alpha}} \right)\\
  & \leqslant & K (\varepsilon, T^{\ast}) e^{- \frac{\nu}{2} t^{\frac{1 -
  \alpha}{1 + \alpha}}} .
\end{eqnarray*}
This completes the proof of
\[ | F (x, t) | \leqslant \max_{D_1} (| g |) + K \frac{\int_{\mathbb{R}
   \backslash D} e^{H (y)} d y}{\int_{\mathbb{R}} e^{H (y)} d y} \leqslant
   \max_D (| g |) + K (\varepsilon, T^{\ast}) e^{- \frac{\nu}{2} t^{\frac{1 -
   \alpha}{1 + \alpha}}} \]
with $D_1 = \left[  \frac{\lambda_1}{2} t^{\frac{1}{1 + \alpha}}, 2 \lambda_2
t^{\frac{1}{1 + \alpha}} \right]$ and $\lambda_1 > 0$ in the case $\frac{| x
|}{t^{\frac{1}{1 + \alpha}}} \leqslant \varepsilon$.
\[ \  \]

\subsubsection{Estimate on $F$ in the case $\frac{x}{t^{\frac{1}{1 + \alpha}}}
\geqslant \varepsilon$}

\

Remark that here, $x \geqslant \varepsilon t^{\frac{1}{1 + \alpha}} \geqslant
2$ for $T^{\ast}$ large enough. From equation (\ref{agh}), we have for any
$y^{\ast} \in \mathbb{R}$ such that $H' (y^{\ast}) = 0$ that $x - y^{\ast}
\leqslant 0$, hence $y^{\ast} \geqslant x$. Furthermore, there exists $\Lambda
> 1$ such that $y^{\ast} \leqslant \Lambda x$. Indeed, if $y^{\ast} \geqslant
\Lambda x$, then equation (\ref{agh}) implies that
\[ - \kappa_2 t \leqslant (1 + | y^{\ast} |)^{\alpha} (x - y^{\ast}) \leqslant
   - (\Lambda - 1) x^{1 + \alpha}, \]
which is in contradiction with $t \leqslant \frac{x^{1 +
\alpha}}{\varepsilon^{1 + \alpha}}$ for $\Lambda$ large enough depending on
$\varepsilon$ (and $T^{\ast}$ large enough).

In summary, in that case, $\Lambda x \geqslant y^{\ast} \geqslant x$ for some
constant $\Lambda > 1$ depending on $\varepsilon$, since $x > 0$ we compute
that for some constant $K > 0$,
\[ H (x) = - \frac{1}{2} \int_0^x f_0 (z) d z \geqslant K x^{1 - \alpha} . \]
This implies that $\max H \geqslant H (x) \geqslant K x^{1 - \alpha} > 0.$

Remark that for $y \leqslant 0$, since $f_0 \leqslant 0$, we have $-
\frac{1}{2} \int_0^y f_0 (z) d z \leqslant 0$, hence
\begin{equation}
  H (y) \leqslant - \frac{(x - y)^2}{4 t} . \label{21778}
\end{equation}
Furthermore, we have in that case $t \leqslant \left( \frac{x}{\varepsilon}
\right)^{1 + \alpha}$, which implies the inequality
\[ \frac{- 1}{t} \leqslant - \left( \frac{\varepsilon}{x} \right)^{1 +
   \alpha}, \]
and we deduce that for $\gamma > 0$,
\begin{eqnarray*}
  H (\gamma x) & \leqslant & \frac{- 1}{4 t} x^2 (1 - \gamma)^2 + K x^{1 -
  \alpha} \gamma^{1 - \alpha}\\
  & \leqslant & x^{1 - \alpha} \left( \frac{- \varepsilon^{1 + \alpha}}{4} (1
  - \gamma)^2 + K \gamma^{1 - \alpha} \right) .
\end{eqnarray*}
Remark that $\frac{- \varepsilon^{1 + \alpha}}{4} (1 - \gamma)^2 + K \gamma^{1
- \alpha} \leqslant 0$ if $\gamma > 0$ is close to $0$ or $\gamma$ is large
(depending on $\varepsilon$). We deduce that there exists $\lambda_2 >
\lambda_1 > 0$ depending on $\varepsilon$ such that for $y$ outside of $D_2
\assign [\lambda_1 x, \lambda_2 x]$,
\begin{equation}
  \frac{H (y)}{\max H} \leqslant \frac{1}{4} \label{319954}
\end{equation}
(we recall that here $\max H \geqslant K x^{1 - \alpha} > 0$). We compute as
in the previous case that
\[ | F (x, t) | \leqslant \max_{D_2} (| g |) + K \frac{\int_{\mathbb{R}
   \backslash D_2} e^{H (y)} d y}{\int_{\mathbb{R}} e^{H (y)} d y} . \]
Here, $x \geqslant \varepsilon t^{\frac{1}{1 + \alpha}}$, therefore $D_2
\subset \left[ \lambda_1 \varepsilon t^{\frac{1}{1 + \alpha}}, + \infty
\right[$ with $\lambda_1 \varepsilon > 0$. Furthermore, we check as in the
previous case that for $z \in [- 1, 1]$ and any $x^{\ast}$ such that $H
(x^{\ast}) = \max H$,
\[ | H (x^{\ast} + z) - H (x^{\ast}) | \leqslant K, \]
therefore
\[ \int_{\mathbb{R}} e^{H (y)} d y \geqslant K e^{\max H} . \]
With (\ref{21778}), we check that
\[ \int_{| y | \leqslant 0} e^{H (y)} d y \leqslant \int_{| y | \leqslant 0}
   e^{- \frac{(x - y)^2}{4 t}} d y \leqslant \int_{\mathbb{R}} e^{-
   \frac{y^2}{4 t}} d y \leqslant K \sqrt{t} . \]
Now, with $H (y) = - \frac{(x - y)^2}{4 t} - \frac{1}{2} \int_0^y f_0 (z) d
z$, we check that there exists $\Lambda_2 > \lambda_2$ depending on
$\varepsilon$ such that if $y \geqslant \Lambda_2 x$, then
\[ H (y) \leqslant \frac{- (x - y)^2}{8 t} . \]
We deduce that
\[ \int_{\Lambda_2 x}^{+ \infty} e^{H (y)} d y \leqslant \int_{\Lambda_2 x}^{+
   \infty} e^{- \frac{(x - y)^2}{8 t}} d y \leqslant \int_{\mathbb{R}} e^{-
   \frac{y^2}{8 t}} d y \leqslant K \sqrt{t} . \]
Finally, with (\ref{319954}) we check that
\[ \int_{[0, \Lambda_2 x] \backslash D_2} e^{H (y)} d y \leqslant
   e^{\frac{1}{4} \max H} \tmop{vol} ([0, \Lambda_2 x] \backslash D_2)
   \leqslant K (\varepsilon, T^{\ast}) x e^{\frac{1}{4} \max H}, \]
where $\tmop{vol} (A)$ is the Lebesgue measure of the set $A$. Combining these
estimates and $\max H \geqslant K x^{1 - \alpha} \geqslant K \varepsilon
t^{\frac{1}{1 + \alpha}}$, we infer that
\begin{eqnarray*}
  \frac{\int_{\mathbb{R} \backslash D_2} e^{H (y)} d y}{\int_{\mathbb{R}} e^{H
  (y)} d y} & \leqslant & K (\varepsilon, T^{\ast}) \frac{\sqrt{t} + x
  e^{\frac{1}{4} \max H}}{e^{\max H}}\\
  & \leqslant & K (\varepsilon, T^{\ast}) \left( \sqrt{t} e^{- \max H} + x
  e^{- \frac{3}{4} \max H} \right)\\
  & \leqslant & K (\varepsilon, T^{\ast})  \left( \sqrt{t} e^{- K \varepsilon
  t^{\frac{1}{1 + \alpha}}} + x e^{- \frac{K}{4} x^{1 - \alpha}} e^{- \frac{K
  \varepsilon}{2} t^{\frac{1}{1 + \alpha}}} \right)\\
  & \leqslant & K (\varepsilon, T^{\ast}) e^{- \frac{K \varepsilon}{2}
  t^{\frac{1}{1 + \alpha}}} .
\end{eqnarray*}
This concludes the proof of
\[ | F (x, t) | \leqslant \max_{D_2} (| g |) + K \frac{\int_{\mathbb{R}
   \backslash D} e^{H (y)} d y}{\int_{\mathbb{R}} e^{H (y)} d y} \leqslant
   \max_D (| g |) + K (\varepsilon, T^{\ast}) e^{- \frac{\nu}{2} t^{\frac{1 -
   \alpha}{1 + \alpha}}} \]
for some small $\nu > 0$ with $D_2 \subset \left[ \lambda_1 \varepsilon
t^{\frac{1}{1 + \alpha}}, + \infty \right[, \lambda_1 \varepsilon > 0$ in the
case $\frac{x}{t^{\frac{1}{1 + \alpha}}} \geqslant \varepsilon$.

\subsubsection{Estimate on $F$ in the case $- \frac{1}{\varepsilon} \leqslant
\frac{x}{t^{\frac{1}{1 + \alpha}}} \leqslant - \varepsilon$}

In that case, for $y^{\ast} \in \mathbb{R}$ solution of $H' (y^{\ast}) = 0$
and $z^{\ast} = \frac{y^{\ast}}{t^{\frac{1}{1 + \alpha}}}$, we recall equation
(\ref{654987}):
\[ - \kappa_2 \leqslant \left( \frac{1}{t^{\frac{1}{1 + \alpha}}} + | z^{\ast}
   | \right)^{\alpha} \left( \frac{x}{t^{\frac{1}{1 + \alpha}}} - z^{\ast}
   \right) \leqslant - \kappa_1 . \]
We deduce that $| z^{\ast} | \geqslant \varepsilon^{\frac{1}{\alpha} + 1}$,
since otherwise $| z^{\ast} | \leqslant \varepsilon^{\frac{1}{\alpha} + 1}$
and
\begin{eqnarray*}
  \kappa_1 & \leqslant & \left| \left( \frac{1}{t^{\frac{1}{1 + \alpha}}} + |
  z^{\ast} | \right)^{\alpha} \left( \frac{x}{t^{\frac{1}{1 + \alpha}}} -
  z^{\ast} \right) \right|\\
  & \leqslant & \left( \frac{1}{t^{\frac{1}{1 + \alpha}}} + \varepsilon^{1 +
  \frac{1}{\alpha}} \right)^{\alpha} \left( \frac{1}{\varepsilon} +
  \varepsilon^{1 + \frac{1}{\alpha}} \right)\\
  & \leqslant & \left( 1 + \frac{1}{t^{\frac{1}{1 + \alpha}} \varepsilon^{1 +
  \frac{1}{\alpha}}} \right) \left( \varepsilon^{\alpha} + \varepsilon^{2 +
  \frac{1}{\alpha} + \alpha} \right)
\end{eqnarray*}
which is impossible if $\varepsilon$ is small enough and $T^{\ast}$ large
enough.

We check also that $| z^{\ast} | \leqslant \frac{2}{\varepsilon}$. Otherwise,
\[ \left| \frac{x}{t^{\frac{1}{1 + \alpha}}} - z^{\ast} \right| \geqslant
   \frac{1}{\varepsilon} \]
and $\left( \frac{1}{t^{\frac{1}{1 + \alpha}}} + | z^{\ast} | \right)^{\alpha}
\geqslant \left( \frac{2}{\varepsilon} \right)^{\alpha},$ thus
\[ \kappa_2 \geqslant \left| \left( \frac{1}{t^{\frac{1}{1 + \alpha}}} + |
   z^{\ast} | \right)^{\alpha} \left( \frac{x}{t^{\frac{1}{1 + \alpha}}} -
   z^{\ast} \right) \right| \geqslant \frac{2^{\alpha}}{\varepsilon^{1 +
   \alpha}}, \]
which is impossible if $\varepsilon$ is small enough. We define as previously
\[ C_0 = t^{- \frac{1 - \alpha}{1 + \alpha}} \max H, \]
and we just showed that for $y^{\ast} \in \mathbb{R}$ solution of $H'
(y^{\ast}) = 0$, we have
\[ | y^{\ast} | \in \left[ \varepsilon^{\frac{1}{\alpha} + 1} t^{\frac{1}{1 +
   \alpha}}, \frac{2}{\varepsilon} t^{\frac{1}{1 + \alpha}} \right] . \]
We can show as previously that $C_0$ is bounded uniformly in time. Now, we
compute that for $\gamma \in \mathbb{R}$,
\[ H' \left( \gamma t^{\frac{1}{1 + \alpha}} \right) \geqslant \frac{t^{-
   \frac{\alpha}{1 + \alpha}}}{2} \left( \frac{x}{t^{\frac{1}{1 + \alpha}}} -
   \gamma + \frac{\kappa_1}{\left( t^{- \frac{1}{1 + \alpha}} + | \gamma |
   \right)^{\alpha}} \right) \]
and
\[ H' \left( \gamma t^{\frac{1}{1 + \alpha}} \right) \leqslant \frac{t^{-
   \frac{\alpha}{1 + \alpha}}}{2} \left( \frac{x}{t^{\frac{1}{1 + \alpha}}} -
   \gamma + \frac{\kappa_2}{\left( t^{- \frac{1}{1 + \alpha}} + | \gamma |
   \right)^{\alpha}} \right) . \]
In particular, there exists $\lambda > 0$ small enough (depending on
$T^{\ast}$ and $\varepsilon$) such that for $\gamma \in [- \lambda, \lambda]$
and $T^{\ast}$ large enough (depending on $\varepsilon$), we have
\begin{equation}
  H' \left( \gamma t^{\frac{1}{1 + \alpha}} \right) \geqslant t^{-
  \frac{\alpha}{1 + \alpha}} . \label{alva}
\end{equation}
This is because $\frac{\kappa_1}{\left( t^{- \frac{1}{1 + \alpha}} + | \gamma
| \right)^{\alpha}} \rightarrow + \infty$ when $| \gamma | \rightarrow 0$ and
$t \rightarrow + \infty$, while $\frac{x}{t^{\frac{1}{1 + \alpha}}} - \gamma$
stays bounded.

In the limit $| \gamma | \rightarrow \infty$, the dominating term is $-
\gamma$, and therefore there exists $\Lambda_1 < 0, \Lambda_2 > 0$ (depending
on $\varepsilon$) such that for all $\gamma < \Lambda_1$,
\begin{equation}
  H' \left( \gamma t^{\frac{1}{1 + \alpha}} \right) \geqslant t^{-
  \frac{\alpha}{1 + \alpha}} \label{140}
\end{equation}
and for all $\gamma > \Lambda_2$,
\begin{equation}
  H' \left( \gamma t^{\frac{1}{1 + \alpha}} \right) \leqslant - t^{-
  \frac{\alpha}{1 + \alpha}} . \label{141}
\end{equation}
We recall that the maximum of $H$ is reached either in $\left[
\varepsilon^{\frac{1}{\alpha} + 1} t^{\frac{1}{1 + \alpha}},
\frac{2}{\varepsilon} t^{\frac{1}{1 + \alpha}} \right]$ or in $\left[ -
\frac{2}{\varepsilon} t^{\frac{1}{1 + \alpha}}, -
\varepsilon^{\frac{1}{\alpha} + 1} t^{\frac{1}{1 + \alpha}} \right]$. Let us
show that (\ref{alva}) implies that for $y \in \left[ - \frac{\lambda}{2}
t^{\frac{1}{1 + \alpha}}, \frac{\lambda}{2} t^{\frac{1}{1 + \alpha}} \right],$
we have
\begin{equation}
  H (y) \leqslant \left( C_0 - \frac{\lambda}{4} \right) t^{\frac{1 -
  \alpha}{1 + \alpha}} \label{123987} .
\end{equation}
Indeed, if it's not true for some $y$, then
\begin{eqnarray*}
  H (\lambda t^{\frac{1}{1+\alpha}}) & = & H (y) + \int_y^{\lambda t^{\frac{1}{1 + \alpha}}} H' (z) d
  z\\
  & \geqslant & \left( C_0 - \frac{\lambda}{4} \right) t^{\frac{1 - \alpha}{1
  + \alpha}} + \int_{\frac{\lambda}{2} t^{\frac{1}{1 + \alpha}}}^{\lambda
  t^{\frac{1}{1 + \alpha}}} t^{- \frac{\alpha}{1 + \alpha}} d z\\
  & \geqslant & \left( C_0 - \frac{\lambda}{4} \right) t^{\frac{1 - \alpha}{1
  + \alpha}} + \frac{\lambda}{2} t^{\frac{1 - \alpha}{1 + \alpha}}\\
  & \geqslant & \left( C_0 + \frac{\lambda}{4} \right) t^{\frac{1 - \alpha}{1
  + \alpha}}\\
  & > & \max H
\end{eqnarray*}
which is a contradiction. We then define the domain
\[ D_3 \assign \left[ (\Lambda_1 - 1) t^{\frac{1}{1 + \alpha}}, (\Lambda_2 +
   1) t^{\frac{1}{1 + \alpha}} \right] \backslash \left[ - \frac{\lambda}{2}
   t^{\frac{1}{1 + \alpha}}, \frac{\lambda}{2} t^{\frac{1}{1 + \alpha}}
   \right] \]
and we estimate as the other steps that
\[ | F (x, t) | \leqslant \left| \frac{\int_{\mathbb{R}} g_0 (y) e^{H (y)} d
   y}{\int_{\mathbb{R}} e^{H (y)} d y} \right| \leqslant \max_{D_3}  | g_0 (y)
   | + K \| f_0 \|_{L^{\infty}} e^{- C_0 t^{\frac{1 - \alpha}{1 + \alpha}}}
   \int_{\mathbb{R} \backslash D_3} e^{H (y)} d y. \]
We have $\mathbb{R} \backslash D_3 = \left[ - \frac{\lambda}{2} t^{\frac{1}{1
+ \alpha}}, \frac{\lambda}{2} t^{\frac{1}{1 + \alpha}} \right] \cup \left(
\mathbb{R} \backslash \left[ (\Lambda_1 - 1) t^{\frac{1}{1 + \alpha}},
(\Lambda_2 + 1) t^{\frac{1}{1 + \alpha}} \right] \right)$ and by
(\ref{123987}),
\[ \int_{\left[ - \frac{\lambda}{2} t^{\frac{1}{1 + \alpha}},
   \frac{\lambda}{2} t^{\frac{1}{1 + \alpha}} \right]} e^{H (y)} d y \leqslant
   e^{\left( C_0 - \frac{\lambda}{4} \right) t^{\frac{1 - \alpha}{1 +
   \alpha}}} \int_{\left[ - \frac{\lambda}{2} t^{\frac{1}{1 + \alpha}},
   \frac{\lambda}{2} t^{\frac{1}{1 + \alpha}} \right]} d y \leqslant \lambda
   t^{\frac{1}{1 + \alpha}} e^{\left( C_0 - \frac{\lambda}{4} \right)
   t^{\frac{1 - \alpha}{1 + \alpha}}} . \]
Furthermore, by equations (\ref{140}) and (\ref{141}), and with similar
argument as the previous steps, we show that for some $\nu > 0$ small and
depending on $\varepsilon, T^{\ast}$,
\[ \int_{\mathbb{R} \backslash \left[ (\Lambda_1 - 1) t^{\frac{1}{1 +
   \alpha}}, (\Lambda_2 + 1) t^{\frac{1}{1 + \alpha}} \right]} e^{H (y)} d y
   \leqslant K (\varepsilon, T^{\ast}) \sqrt{t} e^{(C_0 - \nu) t^{\frac{1 -
   \alpha}{1 + \alpha}}} . \]
We deduce that
\begin{eqnarray*}
  &  & e^{- C_0 t^{\frac{1 - \alpha}{1 + \alpha}}} \int_{\mathbb{R}
  \backslash D_3} e^{H (y)} d y\\
  & \leqslant & K (\varepsilon, T^{\ast}) e^{- \left( \min \left( \nu,
  \frac{\lambda}{2} \right) \right) t^{\frac{1 - \alpha}{1 + \alpha}}} \left(
  1 + \sqrt{t} \right)\\
  & \leqslant & K (\varepsilon, T^{\ast}) e^{- \left( \min \left(
  \frac{\nu}{2}, \frac{\lambda}{4} \right) \right) t^{\frac{1 - \alpha}{1 +
  \alpha}}},
\end{eqnarray*}
which concludes the proof of
\[ | F (x, t) | \leqslant \max_{D_3} (| g |) + K e^{- \left( \min \left(
   \frac{\nu}{2}, \frac{\lambda}{4} \right) \right) t^{\frac{1 - \alpha}{1 +
   \alpha}}} \]
with $D_3 \subset \left] - \infty, - \frac{\lambda}{2} t^{\frac{1}{1 +
\alpha}} \right] \cup \left[ \frac{\lambda}{2} t^{\frac{1}{1 + \alpha}}, +
\infty \right[$ in the case $- \frac{1}{\varepsilon} \leqslant
\frac{x}{t^{\frac{1}{1 + \alpha}}} \leqslant - \varepsilon$.

\subsubsection{Estimate on $F$ in the case $\frac{x}{t^{\frac{1}{1 + \alpha}}}
\leqslant - \frac{1}{\varepsilon}$}

In that case, $x < - 1$ for $T^{\ast}$ large enough and we compute that
\[ H (x) = \frac{- 1}{2} \int_0^x f_0 (z) d z \geqslant - C_1 | x |^{1 -
   \alpha} \]
for some $C_1 > 0$. Now, take $\gamma > 0$, and we estimate, for $T^{\ast}$
large enough, that
\[ - \frac{(1 - \gamma)^2}{4 t} x^2 - K \kappa_2 \gamma^{1 - \alpha} | x |^{1
   - \alpha} \geqslant H (\gamma x) \geqslant - \frac{(1 - \gamma)^2}{4 t} x^2
   - K \kappa_1 \gamma^{1 - \alpha} | x |^{1 - \alpha} . \]
We check that if $\gamma \leqslant \varepsilon$ or $\gamma \geqslant
\frac{1}{\varepsilon}$ for $\varepsilon$ small enough, then
\[ H (\gamma x) \leqslant - 2 C_1 K | x |^{1 - \alpha} \]
(using that $\frac{x}{t^{\frac{1}{1 + \alpha}}} \leqslant -
\frac{1}{\varepsilon}$). We deduce that the maximum of $H$ is reached in $D_4
\assign \left[ \frac{x}{\varepsilon}, \varepsilon x \right]$, and with similar
arguments as the second case, we deduce that for some $\nu > 0$
\[ | F (x, t) | \leqslant \max_{D_4} (| g |) + K (\varepsilon, T^{\ast}) e^{-
   \nu t^{\frac{1 - \alpha}{1 + \alpha}}} \]
with $D_4 \subset \left] - \infty, - \frac{\varepsilon}{2} t^{\frac{1}{1 +
\alpha}} \right]$ in the case $\frac{x}{t^{\frac{1}{1 + \alpha}}} \leqslant -
\frac{1}{\varepsilon}$.

\subsection{Proof of Theorem \ref{mainth} and Proposition \ref{Prop14}}

\begin{proof}[of Theorem \ref{mainth}]
  We recall that the solution of the viscous Burgers equation is
  \[ f (x, t) = \frac{\int_{\mathbb{R}} f_0 (y) e^{H (y)} d
     y}{\int_{\mathbb{R}} e^{H (y)} d y} \]
  where $H (y) = - \frac{(x - y)^2}{4 t} - \frac{1}{2} \int_0^y f_0 (z) d z$.
  The estimate
  \[ \| f (., t) \|_{L^{\infty} (\mathbb{R})} \leqslant \frac{K_2}{(1 +
     t)^{\frac{\alpha}{1 + \alpha}}} \]
  in Theorem \ref{mainth} is a consequence of Proposition \ref{prop211} for $g
  = f_0$ since
  \[ \sup_D | f_0 | \leqslant \frac{K}{(1 + t)^{\frac{\alpha}{1 + \alpha}}} .
  \]

  To show the lower bound, let us compute $f (0, t) < 0$. We have defined in
  the first case (which contains the specific value $x = 0$) of the proof of
  Proposition \ref{prop211} the set $D_1 = \left[  \frac{\lambda_1}{2}
  t^{\frac{1}{1 + \alpha}}, 2 \lambda_2 t^{\frac{1}{1 + \alpha}} \right]$ for
  some $\lambda_1, \lambda_2 > 0$. We estimate
  \begin{eqnarray*}
    | f (0, t) | & \geqslant & \inf_{D_1} (| f_0 |) \frac{\int_{D_1} e^{H (y)}
    d y}{\int_{\mathbb{R}} e^{H (y)} d y} - \| f_0 \|_{L^{\infty}} \left|
    \frac{\int_{\mathbb{R} \backslash D_1} e^{H (y)} d y}{\int_{\mathbb{R}}
    e^{H (y)} d y} \right|\\
    & \geqslant & \inf_{D_1} (| f_0 |) - (\inf_{D_1} (| f_0 |) + \| f_0
    \|_{L^{\infty}}) \left| \frac{\int_{\mathbb{R} \backslash D_1} e^{H (y)} d
    y}{\int_{\mathbb{R}} e^{H (y)} d y} \right|\\
    & \geqslant & \inf_{D_1} (| f_0 |) - 2 \| f_0 \|_{L^{\infty}} \left|
    \frac{\int_{\mathbb{R} \backslash D_1} e^{H (y)} d y}{\int_{\mathbb{R}}
    e^{H (y)} d y} \right|
  \end{eqnarray*}
  and we have shown in the first case that
  \[ \left| \frac{\int_{\mathbb{R} \backslash D_1} e^{H (y)} d
     y}{\int_{\mathbb{R}} e^{H (y)} d y} \right| \leqslant K e^{- \nu
     t^{\frac{1 - \alpha}{1 + \alpha}}} \]
  for some constants $K, \nu > 0$. Since $\inf_{D_1} (| f_0 |) \geqslant
  \frac{K}{t^{\frac{\alpha}{1 + \alpha}}}$, for $t$ large enough,
  \[ | f (0, t) | \geqslant \frac{K_1}{t^{\frac{\alpha}{1 + \alpha}}} . \]
  This completes the proof of Theorem \ref{mainth}.
\end{proof}

\begin{proof}[of Proposition \ref{Prop14}]
  Take a function $g \in C^2 (\mathbb{R}, \mathbb{R})$ with $\| g \|_{C^2}
  \leqslant \kappa$ and define
  \[ A_g (x, t) \assign \int_{\mathbb{R}} g (y) e^{H (y)} d y. \]
  We compute that $\partial_x H = \frac{- 1}{2 t} (x - y)$ and $\frac{- 1}{2
  t} (x - y) e^{- \frac{(x - y)^2}{4 t}} = - \partial_y \left( e^{- \frac{(x -
  y)^2}{4 t}} \right)$, therefore
  \begin{eqnarray*}
    \partial_x A_g & = & - \int_{\mathbb{R}} g (y) \partial_y \left( e^{-
    \frac{(x - y)^2}{4 t}} \right) e^{- \frac{1}{2} \int_0^y f_0 (z) d z} d
    y\\
    & = & \int_{\mathbb{R}} \left( g' (y) - \frac{1}{2} g (y) f_0 (y) \right)
    e^{H (y)} d y\\
    & = & A_{g' - \frac{1}{2} g f_0} .
  \end{eqnarray*}
  Similarly, we have $\partial_t H = \frac{(x - y)^2}{4 t^2}$ and therefore
  \begin{eqnarray*}
    \partial_t A_g & = & \int_{\mathbb{R}} \frac{(x - y)}{2 t} g (y)
    \partial_y \left( e^{- \frac{(x - y)^2}{4 t}} \right) e^{- \frac{1}{2}
    \int_0^y f_0 (z) d z} d y\\
    & = & \int_{\mathbb{R}} \frac{g (y)}{2 t} e^{H (y)} d y\\
    & - & \int_{\mathbb{R}} \left( g' (y) - \frac{1}{2} g (y) f_0 (y) \right)
    \frac{(x - y)}{2 t} e^{H (y)} d y.
  \end{eqnarray*}
  We continue,
  \begin{eqnarray*}
    &  & - \int_{\mathbb{R}} \left( g' (y) - \frac{1}{2} g (y) f_0 (y)
    \right) \frac{(x - y)}{2 t} e^{H (y)} d y\\
    & = & - \int_{\mathbb{R}} \left( g' (y) - \frac{1}{2} g (y) f_0 (y)
    \right) \partial_y \left( e^{- \frac{(x - y)^2}{4 t}} \right) e^{-
    \frac{1}{2} \int_0^y f_0 (z) d z} d y\\
    & = & \int_{\mathbb{R}} \left( g'' (y) - \frac{1}{2} \partial_y (g f_0)
    (y) - \frac{1}{2} g' (y) f_0 (y) + \frac{1}{4} g (y) f_0^2 (y) \right)
    e^{H (y)} d y,
  \end{eqnarray*}
  therefore
  \begin{eqnarray*}
    \partial_t A_g & = & \int_{\mathbb{R}} \left( \frac{g (y)}{2 t} + g'' (y)
    - g' (y) f_0 (y) - \frac{1}{2} g (y) f_0' (y) + \frac{1}{4} g (y) f_0^2
    (y) \right) e^{H (y)} d y\\
    & = & A_{\frac{g}{2 t} + g'' - g' f_0 - \frac{1}{2} g f_0' + \frac{1}{4}
    g f_0^2} .
  \end{eqnarray*}
  Remark that
  \[ f (x, t) = \frac{A_{f_0}}{A_1}, \]
  and therefore $\partial_t^n \partial_x^k f$ can be written as a sum of
  products of terms of the form $\frac{A_g}{A_1}$ for functions $g$ that are
  polynoms on the variables $f_0^{(i)}$ for $i \leqslant 2 n + k$. By
  Proposition \ref{prop211}, we have
  \[ \left| \frac{A_g}{A_1} \right| \leqslant \sup_D | g | + K e^{- \nu
     t^{\frac{1 - \alpha}{1 + \alpha}}} \]
  for $D = \left] - \infty, \mu_1 t^{\frac{1}{1 + \alpha}} \right] \cup \left[
  \mu_2 t^{\frac{1}{1 + \alpha}}, + \infty \right[$ where $\mu_1 < 0, \mu_2 >
  0, \nu, K > 0$ depend on $n, k$. We recall that \ $\partial_x A_g = A_{g' -
  \frac{1}{2} g f_0}$ and $\partial_t A_g = A_{\frac{g}{2 t} + g'' - g' f_0 -
  \frac{1}{2} g f_0' + \frac{1}{4} g f_0^2}$. With the hypothesis on $f_0$ in
  Proposition \ref{Prop14}, we check that for a function $g$ that is
  polynomial on the variables $f_0^{(i)}$ for $i \leqslant 2 n + k - 1$, then
  we will gain in $\sup_D \left| g' - \frac{1}{2} g f_0 \right|$ an additional
  factor $t^{- \frac{\alpha}{1 + \alpha}}$ compared to an estimation on
  $\sup_D | g |$. Similarly, if $i \leqslant 2 n + k - 1$, we will gain in
  $\sup_D \left| \frac{g}{2 t} + g'' - g' f_0 - \frac{1}{2} g f_0' +
  \frac{1}{4} g f_0^2 \right|$ an additional factor $t^{- \frac{2 \alpha}{1 +
  \alpha}}$ compared to $\sup_D | g |$ (since $1 > \frac{2 \alpha}{1 +
  \alpha}$ for $0 < \alpha < 1$).
  
  This concludes the proof of Proposition \ref{Prop14}.
\end{proof}

\section{Computation of the first order}\label{sec4}

This section is devoted to the proofs of Theorem \ref{thfo} and Propositions
\ref{3rdcase}, \ref{lng} and \ref{IDwhynot}. All these proofs are similar in
their approach. We start with the proof of Theorem \ref{thfo} in subsection
\ref{ss32}, and then we will explain what are the differences to get the
propositions in subsection \ref{ss33}.

\subsection{Definition of the functions $y^{\ast}_{\pm}$ and
$y_m^{\ast}$}\label{sec41}

We recall the definition for $\alpha \in] 0, 1 [$ and $\kappa > 0$ of the
function
\[ g (y) = y + \frac{\kappa}{| y |^{\alpha}} . \]
We are interested in the solutions of $z = g (y)$ (see the figure in the
introduction). First, remark that $g (y) \rightarrow + \infty$ when $y
\rightarrow 0^{\pm}, g (y) \rightarrow \pm \infty$ when $y \rightarrow \pm
\infty$ and $g' (y) \rightarrow 1$ when $y \rightarrow \pm \infty$. We compute
for $y \neq 0$ that
\[ g' (y) = 1 - \frac{\kappa \alpha}{y | y |^{\alpha}} . \]
In particular, $g' > 0$ on $] - \infty, 0 [$. We have $g' (y) = 0$ if and only
if $y = y_0 = (\kappa \alpha)^{\frac{1}{1 + \alpha}} .$ This implies that on
$[y_0, + \infty [,$we have $g' (y) > 0$. We compute easily that
\[ g (y_0) = \kappa^{\frac{1}{1 + \alpha}} \left( \alpha^{\frac{1}{1 +
   \alpha}} + \alpha^{- \frac{\alpha}{1 + \alpha}} \right) > 0. \]
By the implicit function theorem, we construct two particular branches of
functions. First, a smooth function $y^{\ast}_- : \mathbb{R} \rightarrow] -
\infty, 0 [$ solution of $z = g (y^{\ast}_- (z))$ for any $z \in \mathbb{R}$,
defined as the inverse of the invertible function $g :] - \infty, 0
[\rightarrow \mathbb{R}$, and another smooth function
\[ y^{\ast}_+ :] g (y_0), + \infty [\rightarrow] y_0, + \infty [ \]
solution of $z = g (y^{\ast}_+ (z))$, defined as the inverse of
\[ g :] y_0, + \infty [\rightarrow] g (y_0), + \infty [. \]
There is also a third branch defined as the inverse of
\[ g :] 0^+, y_0  [\rightarrow] g (y_0), + \infty [ \]
that we denote $y_m^{\ast}$ (it is the one in black in the figure on the
introduction).

\subsection{Proof of Theorem \ref{thfo}}\label{ss32}

\subsubsection{Rescaling in the Hopf-Cole formula}

We recall that the solution of the viscous Burgers equation can be written as
\[ f (x, t) = \frac{\int_{\mathbb{R}} f_0 (y) e^{H (y)} d y}{\int_{\mathbb{R}}
   e^{H (y)} d y} \]
with
\[ H (y) = - \frac{(x - y)^2}{4 t} - \frac{1}{2} \int_0^y f_0 (u) d u. \]
We introduce for $t > 0$ the change of variable $x = z t^{\frac{1}{1 +
\alpha}}$, and we have
\[ f \left( z t^{\frac{1}{1 + \alpha}}, t \right) = \frac{\int_{\mathbb{R}}
   f_0 \left( y t^{\frac{1}{1 + \alpha}} \right) e^{t^{\frac{1 - \alpha}{1 +
   \alpha}} \tilde{H}_t (y, z)} d y}{\int_{\mathbb{R}} e^{t^{\frac{1 -
   \alpha}{1 + \alpha}} \tilde{H}_t (y, z)} d y} \]
where
\[ \tilde{H}_t (y, z) \assign - \frac{(z - y)^2}{4} - \frac{1}{2} t^{- \frac{1
   - \alpha}{1 + \alpha}} \int_0^{y t^{\frac{1}{1 + \alpha}}} f_0 (u) d u. \]
We compute
\[ \partial_y \tilde{H}_t (y, z) = \frac{1}{2} \left( z - y -
   t^{\frac{\alpha}{1 + \alpha}} f_0 \left( y t^{\frac{1}{1 + \alpha}} \right)
   \right) \]
and
\[ \partial_y^2 \tilde{H}_t (y, z) = \frac{1}{2} \left( - 1 - t f_0' \left( y
   t^{\frac{1}{1 + \alpha}} \right) \right) . \]
Remark that $\partial_y^2 \tilde{H}_t$ no longer depends on $z$.

\subsubsection{Construction of the approximate branches $y_{\pm} (., t)$ and
$y_m (., t)$}

We introduce the function
\[ g_t (y) \assign y + t^{\frac{\alpha}{1 + \alpha}} f_0 \left( y
   t^{\frac{1}{1 + \alpha}} \right) . \]
The zeros of $\partial_y \tilde{H}_t$ are the solutions of $z = g_t (y)$. We
compute that
\[ g'_t (y) = 1 + t f_0' \left( y t^{\frac{1}{1 + \alpha}} \right) . \]
First, we remark that for any $\varepsilon > 0$, we have
\[ \| g_t - g \|_{C^1 (\mathbb{R} \backslash [- \varepsilon, \varepsilon])}
   \rightarrow 0 \]
when $t \rightarrow + \infty$, where $g = y + \frac{\kappa}{| y |^{\alpha}}$.
Therefore, we expect some zeros of $\partial_y \tilde{H}_t$ for large time,
i.e. solutions of $z = g_t (y)$, to be close to $y_{\pm}^{\ast} (z)$ or
$y_m^{\ast} (z)$ defined in subsection \ref{sec41}. However, we will not be
able to construct them on the same domains of definition.

\

Take some small $\nu > 0$ (depending on $\alpha$ and $\kappa$ but independent
of $t > 0$). Then, for $t$ large enough, we have that $t f_0' \left( y
t^{\frac{1}{1 + \alpha}} \right) \geqslant 0$ if $y \leqslant - \nu$. In
particular, on $] - \infty, - \nu [$, we have $g_t' \geqslant 1$. We compute
that
\[ g_t (- \nu) = - \nu + \frac{\kappa}{\nu^{\alpha}} (1 + o_{t \rightarrow
   \infty} (1)) \rightarrow + \infty \]
when $\nu \rightarrow 0$, and, for $\nu$ small enough and $t$ large enough, we
have $g_t (- \nu) > 0$. By the implicit function theorem, we deduce that for
$t$ large enough, there exists $y_- (z, t)$ a function solution of $z = g_t
(y_- (z, t))$, which is the inverse of $g_t :] - \infty, - \nu [\rightarrow] -
\infty, g_t (- \nu) [$ defined in
\[ y_- (., t) :] - \infty, g_t (- \nu) [\rightarrow] - \infty, - \nu [. \]
Since $\| g_t - g \|_{C^1 (\mathbb{R} \backslash [- \varepsilon,
\varepsilon])} \rightarrow 0$ for any $\varepsilon > 0$ and $g' (y) > 0$ if $y
> y_0 = (\kappa \alpha)^{\frac{1}{1 + \alpha}} > 0$, for any $\nu > 0$ small
and $t$ large enough, we have that $g'_t (y) > 0$ if $y > y_0 + \nu$. We
therefore construct $y_+ (z, t)$ as the inverse of $g_t :] y_0 + \nu, + \infty
[\rightarrow] g_t (y_0 + \nu), + \infty [$. It is solution of $z = g_t (y_+
(z, t))$ and
\[ y_+ (., t) :] g_t (y_0 + \nu), + \infty [\rightarrow] y_0 + \nu, + \infty
   [. \]
By a similar arguments, if we take $\nu > 0$ small enough and $t$ large enough
(depending on $\nu$), we can construct the middle branch
\[ y_m (., t) : \left] g_t (y_0 + \nu), \frac{1}{\nu} \right[ \rightarrow
   \left] y_0 + \nu, g_t^{- 1} \left( \frac{1}{\nu} \right) \right[ . \]

\subsubsection{Properties of the branches $y_{\pm} (., t)$ and $y_m (., t)$}

We define the set
\[ \mathcal{C}_t \assign \{ (z, y) \in \mathbb{R}^2, z = g_t (y) \}, \]
that is the set of the zeros of $\partial_y \tilde{H}_t$. We want to prove
here that there exists $\delta, \mu, \varepsilon > 0$ small, $T^{\ast} > 0$
large and $C_1, C_2 > 0$, all of them depending only on $\kappa$ and $\alpha$
such that, for $t \geqslant T^{\ast}$,
\begin{enumeratenumeric}
  \item The set $\mathcal{C}_t$ converges to $\mathcal{C}= \{ (z, y) \in
  \mathbb{R}^2, z = g (y) \}$ outside of a vicinity of $\{ (z, y) \in
  \mathbb{R}^2, y = 0 \}$ in the following sense:
  \[ \forall w, \varepsilon > 0, \exists T^{\ast} > 0, \forall (y, z) \in
     \mathcal{C}_t, | y | \geqslant \varepsilon \text{ and } t \geqslant T^{\ast}
     \Rightarrow \exists X \in \mathcal{C}, | X - (y, z) | \leqslant w. \]
  \item
  \[ \| y_- (., t) - y_-^{\ast} \|_{L^{\infty} \left( \left] - \infty,
     \frac{2}{\mu} \right[ \right)} + \| y_+ (., t) - y_+^{\ast}
     \|_{L^{\infty} \left( \left] g (y_0) + \frac{\mu}{2}, + \infty \right[
     \right)} \rightarrow 0 \]
  and
  \[ \| y_m (., t) - y_m^{\ast} \|_{L^{\infty} \left( \left] g (y_0) +
     \frac{\mu}{2}, \frac{2}{\mu} \right[ \right)} \rightarrow 0 \]
  when $t \rightarrow + \infty$.
  
  \item For any $(z, y) \in \left\{ | y | \leqslant t^{- \frac{1}{1 + \alpha}
  + \varepsilon} \right\}$, we have
  \[ | \partial_y \tilde{H}_t (y, z) | \leqslant 1 + | z | + \| f_0
     \|_{L^{\infty}} t^{\frac{\alpha}{1 + \alpha}} . \]
  \item Outside of $\left\{ | y | \leqslant t^{- \frac{1}{1 + \alpha} +
  \varepsilon} \right\}$ we have that $\partial_y \tilde{H}_t (y, z)
  \rightarrow \frac{1}{2} (z - g (y))$ locally uniformly.
  
  \item The set $\mathcal{C}_t \assign \{ (z, y) \in \mathbb{R}^2, \partial_y
  \tilde{H}_t (y, z) = 0 \}$ cut the plane in two parts. In the one containing
  $(0, 1)$, we have $\partial_y \tilde{H}_t < 0$ and in the other one, we have
  $\partial_y \tilde{H}_t > 0$.
  
  \item The set $\mathcal{C}_t \cap \{ z \leqslant g (y_0) + \mu \}$ contains
  parts of the branch $y_-$, and any other elements must be in $\left\{ | y |
  \leqslant t^{- \frac{1}{1 + \alpha} + \varepsilon} \right\}$ or in $B (A,
  \varepsilon)$ where $A = (g (y_0), y_0)$.
  
  \item The set $\mathcal{C}_t \cap \left\{ z \geqslant \frac{1}{\mu}
  \right\}$ contains parts of the branch $y_+$, and any other elements must be
  in $\{ | y | \leqslant 1 \}$.
  
  \item The set $\mathcal{C}_t \cap \left\{ g (y_0) + \mu \leqslant z
  \leqslant \frac{1}{\mu} \right\}$ contains parts of the branches $y_+, y_-,
  y_m$, and any other elements must be in $\left\{ | y | \leqslant t^{-
  \frac{1}{1 + \alpha} + \varepsilon} \right\}$.
  
  \item For $z \in \left[ g (y_0) + \mu, \frac{1}{\mu} \right]$ and $y \in
  \left[ \frac{- 1}{\delta}, - \delta \right] \cup \left[ y_0 + \delta,
  \frac{1}{\delta} \right]$, then
  \[ - C_2 \leqslant \partial_y^2 \tilde{H}_t (y, z) \leqslant - C_1 . \]
\end{enumeratenumeric}
These properties can be summarized in the following way:

\begin{figure}[H]
    \centering
    \includegraphics[width=20cm]{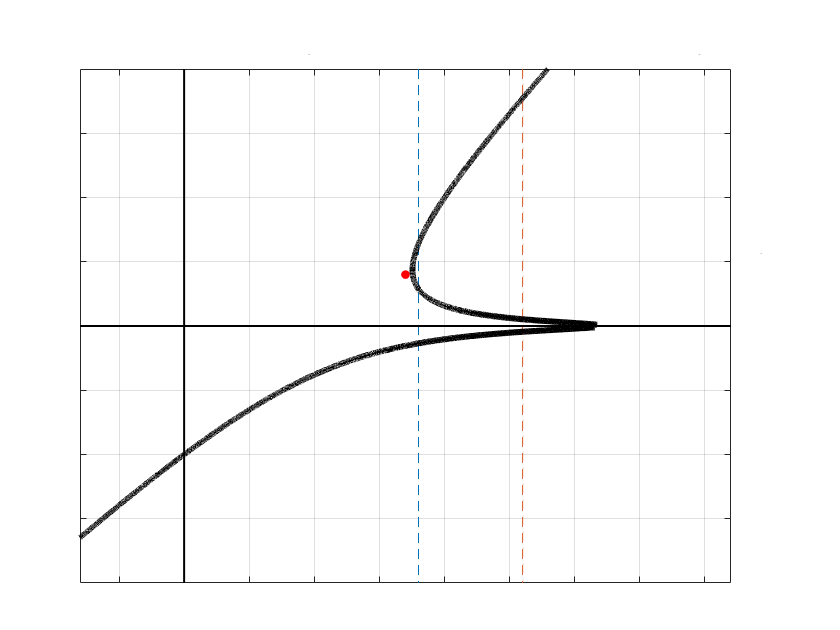}
\end{figure}

\begin{center}
  Plot in the plane $(z, y) \in \mathbb{R}^2$ of $\mathcal{C}_t$ for $\kappa =
  1, \alpha = \frac{1}{3}, t = 0.03$. The red dot is $A = (g (y_0), y_0)$, the
  dotted blue line is $\{ z = g (y_0) + \mu \}$ and the dotted red line is
  $\left\{ z = \frac{1}{\mu} \right\}$. 
\end{center}

\

Let us explain what these properties means on this graph. The set
$\mathcal{C}_t$, that is the set of the zeros of $\partial_y \tilde{H}_t$, in
black on the graph, converges in a sense to the first figure in subsection
\ref{s111} when $t \rightarrow + \infty$, but not uniformly near $\{ y = 0
\}$. In fact, $\mathcal{C}_t$ has only one branch while $\mathcal{C}$ is two
separated branches (we can check that $\mathcal{C}_{} \cap \{ y = 0 \} =
\emptyset$). However, if we avoid $\{ | y | \leqslant \varepsilon \}$ for any
$\varepsilon > 0$, then the convergence will be uniform (see property 1.).

\

As such, the three branches $y_+, y_-, y_m$, that are parts of
$\mathcal{C}_t$ and that can be defined outside of $\{ | y | \leqslant
\varepsilon \}$ for any $\varepsilon > 0$ and outside of a neighborhood of
$A$, the red point, converges to their respective limits $y_+^{\ast},
y_-^{\ast}, y^{\ast}_m$ outside of this domain (see property 2.). However, the
set $\mathcal{C}_t$ can be very complicated near $\{ y = 0 \}$ and near $A$.
Its description there does not only depend on the equivalents at $\pm \infty$
of $f_0$, but on the full function itself. Luckily, we will not need a good
description of $\mathcal{C}_t$ there.

\

More generally than its set of zero, the function $\partial_y \tilde{H}_t$
converges for large time to its limit $\frac{1}{2} (z - g (y))$ uniformly,
except in a shrinking neighborhood of $\{ y = 0 \}$. In it, we can still bound
its values for large time (see properties 3.,4.). We can also compute its sign
on both sides of the plane cut by $\mathcal{C}_t$ (see property 5.)

\

Now, in $\left\{ z \leqslant g(y_0) + \mu \right\}$,
that is left of the dotted blue line, we see that the set $\mathcal{C}_t$
contains the branch $y_-$ there, but also possibly elements near the red point
$A$ or near $\{ y = 0 \}$ (see property 6.). Indeed, although $\mathcal{C}_t$
converge to $\mathcal{C}$ near $A$, it is not obvious that there is only one
smooth cusp there.

\

In the set $\left\{ z \geqslant \frac{1}{\mu} \right\}$, that is right of the
dotted red line, the curve $\mathcal{C}_t$ has not yet fully converges to
$\mathcal{C}$ near $\{ y = 0 \}$. However, it contains the branch $y_+$, which
is far away from $\{ y = 0 \}$ (at distance of size $\frac{1}{\mu}$ for small
$\mu$), and any other element must be in a neighborhood of size $1$ around $\{
y = 0 \}$ (see property 7.).

\

Finally, in the set $\left\{ g (y_0) + \mu \leqslant z \leqslant
\frac{1}{\mu} \right\}$, that is between the two dotted lines, the set
$\mathcal{C}_t$ contains the three branches $y_+, y_-, y_m$ that are all three
uniformly far away from $\{ y = 0 \}$, and any other elements must be in a
shrinking neighborhood of $\{ y = 0 \}$ (see property 8.). We conclude with
property 9., stating that in a vicinity of the branches $y_+, y_-$ and $y_m$
between the two dotted lines, the quantity $\partial_y^2 \tilde{H}_t$ is
uniformly strictly negative.

\

Let us now prove all these properties.

\

Properties 1. and 2. are a direct consequence of $\| g_t - g \|_{C^1
(\mathbb{R} \backslash [- \varepsilon, \varepsilon])} \rightarrow 0$ for any
$\varepsilon > 0$ and the fact that $| g_t' | + | g' | \leqslant K
(\varepsilon)$ on $\mathbb{R} \backslash [- \varepsilon, \varepsilon]$.
Property 3. follows from the computation
\[ | \partial_y \tilde{H}_t (y, z) | = \left| \frac{1}{2} \left( z - y -
   t^{\frac{\alpha}{1 + \alpha}} f_0 \left( y t^{\frac{1}{1 + \alpha}} \right)
   \right) \right| \leqslant \left( | y | + | z | + t^{\frac{\alpha}{1 +
   \alpha}} \| f_0 \|_{L^{\infty} (\mathbb{R})} \right) . \]
Property 4. is a consequence of the fact that if $| y | \geqslant t^{-
\frac{1}{1 + \alpha} + \varepsilon}$, then $\left| y t^{\frac{1}{1 + \alpha}}
\right|$ is large when $t \rightarrow + \infty$, and thus $f_0 \left( y
t^{\frac{1}{1 + \alpha}} \right)$ behaves like its equivalent when the time
becomes large.

\

Property 5. follows from the continuity of $g_t$ and the computations
\[ \partial_y \tilde{H}_t (0, 1) = \frac{- 1 - t^{\frac{\alpha}{1 + \alpha}}
   f_0 \left( t^{\frac{1}{1 + \alpha}} \right)}{2} < 0 \]
for $t$ large enough or $\partial_y \tilde{H}_t (z, 0) \rightarrow - \infty$
when $z \rightarrow - \infty$.

\

Properties 6.,7. and 8. are then a consequence of properties 1.,2. and 4. For
the last property, we recall that
\[ \partial_y^2 \tilde{H}_t (y, z) = \frac{1}{2} \left( - 1 - t f_0' \left( y
   t^{\frac{1}{1 + \alpha}} \right) \right) . \]
In $\{ | y | \geqslant \delta \}$, we have $\frac{1}{2} \left( - 1 - t f_0'
\left( y t^{\frac{1}{1 + \alpha}} \right) \right) \rightarrow \frac{1}{2}
\left( - 1 + \frac{\kappa \alpha}{y | y |^{\alpha}} \right)$ uniformly when $t
\rightarrow + \infty$. We check easily that $- 1 + \frac{\kappa \alpha}{y | y
|^{\alpha}} = 0$ only happens if $y = y_0$, and that $- 1 + \frac{\kappa
\alpha}{y | y |^{\alpha}} < K (\delta) < 0$ if $y \in \left[ y_0 + \delta,
\frac{1}{\delta} \right]$.

\

This completes the proof of the nine properties.

\subsubsection{Position of the maximum of $\tilde{H}_t$}\label{ss324}

Take $z \leqslant g (y_0) + \mu$. Then, the maximum of $y \rightarrow
\tilde{H}_t (y, z)$ must be reached at a zero of $\partial_y \tilde{H}_t (.,
z)$. By 6., this is either at $y = y_- (z, t)$, or possibly in $\left\{ | y |
\leqslant t^{- \frac{1}{1 + \alpha} + \varepsilon} \right\}$ or in $B (A,
\varepsilon)$. Now, by 3. and 4., remark that $y \rightarrow \partial_y
\tilde{H}_t (., z)$ is
\begin{itemizedot}
  \item Strictly increasing on $] - \infty, y_- (z, t) [$
  
  \item Strictly decreasing on $\left] y_- (z, t), - t^{- \frac{1}{1 + \alpha}
  + \varepsilon} \right[$
  
  \item Bounded by $1 + | z | + \| f_0 \|_{L^{\infty}} t^{\frac{\alpha}{1 +
  \alpha}}$ on $\left] - t^{- \frac{1}{1 + \alpha} + \varepsilon}, t^{-
  \frac{1}{1 + \alpha} + \varepsilon} \right[$
  
  \item Strictly decreasing on $\left] t^{- \frac{1}{1 + \alpha} +
  \varepsilon}, + \infty \right[$ except eventually on an interval of size
  $\varepsilon$ near $A$, where it is bounded by a universal constant,
  independent of $\varepsilon$ or $t$.
\end{itemizedot}
Still by 3., there exists $\nu > 0$ depending only on $\mu, \kappa, \alpha$
such that
\[ \tilde{H}_t (y_- (z, t), z) - \tilde{H}_t \left( - t^{- \frac{1}{1 +
   \alpha} + \varepsilon}, z \right) \geqslant \nu \]
and we check easily that
\[ \left| \tilde{H}_t \left( - t^{- \frac{1}{1 + \alpha} + \varepsilon}, z
   \right) - \tilde{H}_t \left( t^{- \frac{1}{1 + \alpha} + \varepsilon}, z
   \right) \right| \leqslant 2 t^{- \varepsilon} \| f_0 \|_{L^{\infty}} + 2 (1
   + | z |) t^{- \frac{1}{1 + \alpha} + \varepsilon} \rightarrow 0 \]
when $t \rightarrow + \infty$.

We deduce that for $z \leqslant g (y_0) + \mu$, the maximum of $y \rightarrow
\tilde{H}_t (y, z)$ is reached at $y_- (z, t)$, and that there exists a
constant $\nu_0 > 0$ depending only on $\alpha, \kappa$ such that
\[ \tilde{H}_t (y_- (z, t), z) \geqslant \tilde{H}_t (y, z) + \nu_0 \]
for any $y \in \left\{ | y | \leqslant t^{- \frac{1}{1 + \alpha} +
\varepsilon} \right\} \cup B (A, \varepsilon)$.

Indeed, the maximum cannot be near $\{ y = 0 \}$ because there $\tilde{H}_t$
cannot increase by $\nu > 0$, and it cannot be near $A$ if we take
$\varepsilon > 0$ small enough.

\

A similar argument can be made for $z \geqslant \frac{1}{\mu}$ if we take
$\mu > 0$ large enough. In that case, the maximum is reached either at $y_+
(z, t)$ or in $\{ | y | \leqslant 1 \}$, but using the above properties we can
show that in $\{ | y | \leqslant 1 \}$,
\[ \tilde{H}_t (y, z) \leqslant - \frac{(z + 1)^2}{4} - K \]
for some universal constant $K > 0$. Now, we estimate
\[ \tilde{H}_t (y_+ (z, t), z) \geqslant - \frac{(z - y_+ (z, t))^2}{4} - K |
   y_+ (z, t) |^{1 - \alpha}, \]
and since $| y_+ (z, t) - z | \rightarrow 0$ when $z \rightarrow + \infty$
(uniformly in $t$), for $z \geqslant \frac{1}{\mu}$ and $\mu$ small enough we
have
\[ \tilde{H}_t (y_+ (z, t), z) \geqslant \tilde{H}_t (y, z) + \nu_0 \]
for a universal constant $\nu_0 > 0$, for any $y \in \mathbb{R}$ such that $|
y | \leqslant 1.$

\

We now focus on the case $g (y_0) + \mu \leqslant z \leqslant \frac{1}{\mu}$.
By similar arguments on the sign of $\partial_y \tilde{H}_t (., z)$ in
different regions, we can show that the maximum is reached at either $y_+ (z,
t)$ or $y_- (z, t)$, with a margin $\nu_0 > 0$ compared to $y_m (z, t)$ (it is
a local minima of $\tilde{H}_t$) and $\left\{ | y | \leqslant t^{- \frac{1}{1
+ \alpha} + \varepsilon} \right\}$.

We introduce the quantities
\[ \mathfrak{h}_+ (z) \assign \tilde{H}_t (y_+ (z, t), z), \mathfrak{h}_- (z)
   \assign \tilde{H}_t (y_- (z, t), z) . \]
Since $\partial_y \tilde{H} (y_{\pm} (z, t), z) = 0$ by construction, we have
\[ \mathfrak{h}_{\pm}' (z) = \partial_z \tilde{H} (y_{\pm} (z, t), z) = -
   \frac{1}{2} (z - y_{\pm} (z, t)) = - t^{\frac{\alpha}{1 + \alpha}} f_0
   \left( y_{\pm} (z, t) t^{\frac{1}{1 + \alpha}} \right) < 0. \]
Also,
\[ \mathfrak{h}_-' (z) -\mathfrak{h}_+' (z) = \frac{1}{2} (- y_+ (z, t) + y_-
   (z, t)) > 0 \]
since $y_- (z, t) < 0$ and $y_+ (z, t) > 0$. We deduce that $\mathfrak{h}_+'
(z) <\mathfrak{h}_-' (z)$ on $\left[ g (y_0) + \mu, \frac{1}{\mu} \right]$.

Since $\mathfrak{h}_+ (1 / \mu) >\mathfrak{h}_- (1 / \mu)$ and $\mathfrak{h}_+
(g (y_0) + \mu) <\mathfrak{h}_- (g (y_0) + \mu)$, by continuity there exists
$z_c (t) \in \left[ g (y_0) + \mu, \frac{1}{\mu} \right]$ such that
$\mathfrak{h}_+ (z_c (t)) =\mathfrak{h}_- (z_c (t))$. That is, $\max
\tilde{H}_t (., z_c (t))$ is reached at two points, $y_+ (z_c (t), t) > 0$ and
$y_- (z_c (t), t) < 0.$ But since $\mathfrak{h}_+' (z) <\mathfrak{h}_-' (z)$
on $\left[ g (y_0) + \mu, \frac{1}{\mu} \right]$, there is only one such point
$z_c (t)$.

By the uniform convergence $\| g_t - g \|_{C^1 (\mathbb{R} \backslash [-
\varepsilon, \varepsilon])} \rightarrow 0$ when $t \rightarrow + \infty$ (the
branches $y_+ (z, t)$ and $y_- (z, t)$ are uniformly far away from $0$ for $z
\in \left[ g (y_0) + \mu, \frac{1}{\mu} \right]$), the quantity $z_c (t)$
converges to a limit $z_c$ when $t \rightarrow + \infty$, define as the unique
solution of $\mathcal{H} (y_+^{\ast} (z_c)) =\mathcal{H} (y_-^{\ast} (z_c))$
where
\[ \mathcal{H} (y) = - \frac{\kappa^2}{4 | y |^{2 \alpha}} - \frac{\kappa (1 -
   \alpha)}{2} | y |^{1 - \alpha} . \]
Indeed, the equation satisfied by $z_c (t)$ is $\tilde{H}_t (y_+ (z_c (t), t),
z_c (t)) = \tilde{H}_t (y_- (z_c (t), t), z_c (t))$, and we get the limit
equation by taking $t \rightarrow + \infty$ using
\[ y_{\pm} (z_c (t), t) - g_t (y_{\pm} (z_c (t), t)) = - t^{\frac{\alpha}{1 +
   \alpha}} f_0 \left( y_{\pm} (z_c (t), t) t^{\frac{1}{1 + \alpha}} \right)
\]
and the fact that $y_{\pm} (z_c (t), t)$ are far from $0$ uniformly in time.

\

From now on, we consider any small $\varepsilon > 0$ and $z \in \mathbb{R}
\backslash [z_c - \varepsilon, z_c + \varepsilon]$. Then, for $t \geqslant
T^{\ast}$ large enough (depending on $\varepsilon$), there exists $\nu
(\varepsilon) > 0$ a small constant such that
\begin{itemizedot}
  \item if $z \leqslant z_c - \varepsilon$, then $\max \tilde{H}_t (., z)$ is
  reached only at $y_- (z, t)$, and $\max \tilde{H}_t (., t) = \tilde{H}_t
  (y_- (z, t), t) \geqslant \tilde{H}_t (y_+ (z, t), t) + \nu (\varepsilon)$
  
  \item if $z \geqslant z_c + \varepsilon,$then $\max \tilde{H}_t (., z)$ is
  reached only at $y_+ (z, t)$, and $\max \tilde{H}_t (., t) = \tilde{H}_t
  (y_+ (z, t), t) \geqslant \tilde{H}_t (y_- (z, t), t) + \nu (\varepsilon)$.
\end{itemizedot}
To simplify the notations, we define for $z \in \mathbb{R} \backslash \{ z_c
\}$
\[ y_{\max} (z, t) \assign \left\{\begin{array}{l}
     y_+ (z, t) \tmop{if} z > z_c\\
     y_- (z, t) \tmop{if} z < z_c
   \end{array}\right. \]
as well as $y_{\max}^{\ast} (z)$, its limit when $t \rightarrow + \infty$, 
and
\[ \tilde{H}_{\max} (z, t) \assign \tilde{H}_t (y_{\max} (z, t), t) . \]
Remark that $| y_{\max} (z, t) | \geqslant C_0 > 0$ where $C_0$ is a constant
independent of $z$ and $t$.

\

Now, we infer that there exists small constants $\gamma, \nu > 0$ depending
on $\alpha, \kappa, \varepsilon$ such that for $t \geqslant T^{\ast}$ large
enough, if $y \in \mathbb{R} \backslash [y_{\max} (z, t) - \gamma, y_{\max}
(z, t) + \gamma]$ and $z \in \mathbb{R} \backslash [z_c - \varepsilon, z_c +
\varepsilon]$, then
\[ \tilde{H}_{\max} (z, t) \geqslant \tilde{H}_t (y, t) + \nu . \]
This is a consequence of 3.,4. and the comparison between $\tilde{H}_t (y_-
(z, t), t)$ and $\tilde{H}_t (y_+ (z, t), t)$ when $g (y_0) + \mu \leqslant z
\leqslant \frac{1}{\mu}$.

\

\subsubsection{End of the proof of Theorem \ref{thfo}}\label{ss325}

We are now equiped to estimate
\[ f \left( z t^{\frac{1}{1 + \alpha}}, t \right) = \frac{\int_{\mathbb{R}}
   f_0 \left( y t^{\frac{1}{1 + \alpha}} \right) e^{t^{\frac{1 - \alpha}{1 +
   \alpha}} \tilde{H}_t (y, z)} d y}{\int_{\mathbb{R}} e^{t^{\frac{1 -
   \alpha}{1 + \alpha}} \tilde{H}_t (y, z)} d y} . \]
We define
\[ \mathcal{A}_{\mathfrak{f}} (z, t) \assign \int_{\mathbb{R}} \mathfrak{f}
   \left( y t^{\frac{1}{1 + \alpha}} \right) e^{t^{\frac{1 - \alpha}{1 +
   \alpha}} \tilde{H}_t (y, z)} d y \]
for $\mathfrak{f} \in \{ 1, f_0 \}$, as well as
\[ \mathcal{D}_{z, t, \gamma} \assign [y_{\max} (z, t) - \gamma, y_{\max} (z,
   t) + \gamma] . \]
We decompose
\[ \mathcal{A}_{\mathfrak{f}} (z, t) = \int_{\mathcal{D}_{z, t, \gamma}}
   \mathfrak{f} \left( y t^{\frac{1}{1 + \alpha}} \right) e^{t^{\frac{1 -
   \alpha}{1 + \alpha}} \tilde{H}_t (y, z)} d y + \int_{\mathbb{R} \backslash
   \mathcal{D}_{z, t, \gamma}} \mathfrak{f} \left( y t^{\frac{1}{1 + \alpha}}
   \right) e^{t^{\frac{1 - \alpha}{1 + \alpha}} \tilde{H}_t (y, z)} d y. \]
On $\mathbb{R} \backslash \mathcal{D}_{z, t, \gamma}$ we have shown that
$\tilde{H}_t (., z) \leqslant \tilde{H}_{\max} (z, t) - \nu$, and using
arguments similar as for the proof of (\ref{123}) to deal with the
integrability, we check that
\begin{eqnarray}
  &  & \left| \int_{\mathbb{R} \backslash \mathcal{D}_{z, t, \gamma}}
  \mathfrak{f} \left( y t^{\frac{1}{1 + \alpha}} \right) e^{t^{\frac{1 -
  \alpha}{1 + \alpha}} \tilde{H}_t (y, z)} d y \right| \nonumber\\
  & \leqslant & K \| \mathfrak{f} \|_{L^{\infty} (\mathbb{R})} e^{t^{\frac{1
  - \alpha}{1 + \alpha}} \left( \tilde{H}_{\max} (z, t) - \frac{\nu}{2}
  \right)} \nonumber\\
  & \leqslant & K e^{t^{\frac{1 - \alpha}{1 + \alpha}} \left(
  \tilde{H}_{\max} (z, t) - \frac{\nu}{2} \right)} .  \label{jun2}
\end{eqnarray}
We continue with
\begin{eqnarray*}
  &  & \int_{\mathcal{D}_{z, t, \gamma}} f_0 \left( y t^{\frac{1}{1 +
  \alpha}} \right) e^{t^{\frac{1 - \alpha}{1 + \alpha}} \tilde{H}_t (y, z)} d
  y\\
  & = & f_0 \left( y_{\max} (z, t) t^{\frac{1}{1 + \alpha}} \right)
  \int_{\mathcal{D}_{z, t, \gamma}} e^{t^{\frac{1 - \alpha}{1 + \alpha}}
  \tilde{H}_t (y, z)} d y\\
  & + & \int_{\mathcal{D}_{z, t, \gamma}} \left( f_0 \left( y t^{\frac{1}{1 +
  \alpha}} \right) - f_0 \left( y_{\max} (z, t) t^{\frac{1}{1 + \alpha}}
  \right) \right) e^{t^{\frac{1 - \alpha}{1 + \alpha}} \tilde{H}_t (y, z)} d
  y.
\end{eqnarray*}
From 9., we have for $y \in \mathcal{D}_{z, t, \gamma}$ (and $\gamma$ small
enough) that
\[ \tilde{H}_t (y, z) \geqslant \tilde{H}_{\max} (z, t) - C_2 (y - y_{\max}
   (z, t))^2, \]
hence
\begin{eqnarray}
  &  & \int_{\mathcal{D}_{z, t, \gamma}} e^{t^{\frac{1 - \alpha}{1 + \alpha}}
  \tilde{H}_t (y, z)} d y \nonumber\\
  & \geqslant & e^{\tilde{H}_{\max} (z, t) t^{\frac{1 - \alpha}{1 + \alpha}}}
  \int_{\mathcal{D}_{z, t, \gamma}} e^{- C_2 t^{\frac{1 - \alpha}{1 + \alpha}}
  (y - y_{\max} (z, t))^2} d y \nonumber\\
  & \geqslant & e^{\tilde{H}_{\max} (z, t) t^{\frac{1 - \alpha}{1 + \alpha}}}
  \int_{- \gamma}^{\gamma} e^{- C_2 t^{\frac{1 - \alpha}{1 + \alpha}} y^2} d y
  \nonumber\\
  & \geqslant & K t^{- \frac{1 - \alpha}{2 (1 + \alpha)}} e^{\tilde{H}_{\max}
  (z, t) t^{\frac{1 - \alpha}{1 + \alpha}}} .  \label{jun1}
\end{eqnarray}
We continue. Still from 9., we have for $y \in \mathcal{D}_{z, t, \gamma}$
(and $\gamma$ small enough) that
\[ \tilde{H}_t (y, z) \leqslant \tilde{H}_{\max} (z, t) - C_1 (y - y_{\max}
   (z, t))^2, \]
hence
\begin{eqnarray}
  &  & \left| \int_{\mathcal{D}_{z, t, \gamma}} \left( f_0 \left( y
  t^{\frac{1}{1 + \alpha}} \right) - f_0 \left( y_{\max} (z, t) t^{\frac{1}{1
  + \alpha}} \right) \right) e^{t^{\frac{1 - \alpha}{1 + \alpha}} \tilde{H}_t
  (y, z)} d y \right| \nonumber\\
  & \leqslant & e^{t^{\frac{1 - \alpha}{1 + \alpha}} \tilde{H}_{\max} (z, t)}
  \int_{\mathcal{D}_{z, t, \gamma}} \left| f_0 \left( y t^{\frac{1}{1 +
  \alpha}} \right) - f_0 \left( y_{\max} (z, t) t^{\frac{1}{1 + \alpha}}
  \right) \right| e^{- C_1 t^{\frac{1 - \alpha}{1 + \alpha}} (y - y_{\max} (z,
  t))^2} d y \nonumber\\
  & \leqslant & e^{t^{\frac{1 - \alpha}{1 + \alpha}} \tilde{H}_{\max} (z, t)}
  \sup_{\mathcal{D}_{z, t, \gamma}} \left| f_0' \left( y t^{\frac{1}{1 +
  \alpha}} \right) \right| \int_{\mathcal{D}_{z, t, \gamma}} t^{\frac{1}{1 +
  \alpha}} | y - y_{\max} | e^{- C_1 t^{\frac{1 - \alpha}{1 + \alpha}} (y -
  y_{\max} (z, t))^2} d y \nonumber\\
  & \leqslant & K e^{t^{\frac{1 - \alpha}{1 + \alpha}} \tilde{H}_{\max} (z,
  t)} t^{- 1 + \frac{1}{1 + \alpha}} \int_{- \gamma}^{\gamma} | y | e^{- C_1
  t^{\frac{1 - \alpha}{1 + \alpha}} y^2} d y \nonumber\\
  & \leqslant & K e^{t^{\frac{1 - \alpha}{1 + \alpha}} \tilde{H}_{\max} (z,
  t)} t^{- \frac{1}{1 + \alpha}} .  \label{jun3}
\end{eqnarray}
We define
\[ \mathcal{I} (z, t) \assign \int_{\mathcal{D}_{z, t, \gamma}} e^{t^{\frac{1
   - \alpha}{1 + \alpha}} \tilde{H}_t (y, z)} d y \]
and $\mathcal{J} (z, t), \mathcal{K} (z, t)$ by
\[ \int_{\mathbb{R}} f_0 \left( y t^{\frac{1}{1 + \alpha}} \right)
   e^{t^{\frac{1 - \alpha}{1 + \alpha}} \tilde{H}_t (y, z)} d y = f_0 \left(
   y_{\max} (z, t) t^{\frac{1}{1 + \alpha}} \right) \mathcal{I} (z, t)
   +\mathcal{J} (z, t) \]
and
\[ \mathcal{K} (z, t) \assign \int_{\mathbb{R}} e^{t^{\frac{1 - \alpha}{1 +
   \alpha}} \tilde{H}_t (y, z)} d y -\mathcal{I} (z, t) . \]
We have
\[ f \left( z t^{\frac{1}{1 + \alpha}}, t \right) = f_0 \left( y_{\max} (z, t)
   t^{\frac{1}{1 + \alpha}} \right) + \frac{\mathcal{J} (z, t) - f_0 \left(
   y_{\max} (z, t) t^{\frac{1}{1 + \alpha}} \right) \mathcal{K} (z,
   t)}{\mathcal{I} (z, t) +\mathcal{K} (z, t)} . \]
From (\ref{jun1}) we have
\[ \mathcal{I} (z, t) \geqslant K t^{- \frac{1 - \alpha}{2 (1 + \alpha)}}
   e^{\tilde{H}_{\max} (z, t) t^{\frac{1 - \alpha}{1 + \alpha}}} \]
and from (\ref{jun2}) we have
\[ | \mathcal{K} (z, t) | \leqslant K e^{- \frac{\nu}{2} t^{\frac{1 -
   \alpha}{1 + \alpha}}} e^{t^{\frac{1 - \alpha}{1 + \alpha}} \tilde{H}_{\max}
   (z, t)} . \]
Therefore, for $t$ large enough, $\mathcal{I} (z, t) +\mathcal{K} (z, t)
\geqslant \mathcal{I} (z, t) / 2$. From (\ref{jun2}) and (\ref{jun3}), for $t$
large enough we have
\[ | \mathcal{J} (z, t) | \leqslant K e^{t^{\frac{1 - \alpha}{1 + \alpha}}
   \tilde{H}_{\max} (z, t)} t^{- \frac{1}{1 + \alpha}} \]
hence
\[ \frac{| \mathcal{J} (z, t) |}{| \mathcal{I} (z, t) +\mathcal{K} (z, t) |}
   \leqslant K \frac{t^{- \frac{1}{1 + \alpha}}}{t^{- \frac{1 - \alpha}{2 (1 +
   \alpha)}}} \leqslant K t^{- \frac{1}{2}} . \]
We check similarly that
\[ \frac{\left| f_0 \left( y_{\max} (z, t) t^{\frac{1}{1 + \alpha}} \right)
   \mathcal{K} (z, t) \right|}{| \mathcal{I} (z, t) +\mathcal{K} (z, t) |}
   \leqslant K e^{- \frac{\nu}{4} t^{\frac{1 - \alpha}{1 + \alpha}}} . \]
Combining these estimates, we have
\[ \left| t^{\frac{\alpha}{1 + \alpha}} f \left( z t^{\frac{1}{1 + \alpha}}, t
   \right) - t^{\frac{\alpha}{1 + \alpha}} f_0 \left( y_{\max} (z, t)
   t^{\frac{1}{1 + \alpha}} \right) \right| \leqslant K t^{- \frac{1}{2} +
   \frac{\alpha}{1 + \alpha}} . \]
Since $| y_{\max} (z, t) | \geqslant C_0 > 0$ where $C_0$ depends only on
$\alpha$ and $\kappa$, we have
\[ \lim_{t \rightarrow + \infty} t^{\frac{\alpha}{1 + \alpha}} f_0 \left(
   y_{\max} (z, t) t^{\frac{1}{1 + \alpha}} \right) = \kappa | y_{\max}^{\ast}
   (z) |^{- \alpha} . \]
This concludes the proof of Theorem \ref{thfo}.

\

About the case $z = z_c$. The difficulty is now that
\[ \int_{y_+ (z, t) - \gamma}^{y_+ (z, t) + \gamma} \mathfrak{f} \left( y
   t^{\frac{1}{1 + \alpha}} \right) e^{t^{\frac{1 - \alpha}{1 + \alpha}}
   \tilde{H}_t (y, z)} d y \tmop{and} \int_{y_- (z, t) - \gamma}^{y_- (z, t) +
   \gamma} \mathfrak{f} \left( y t^{\frac{1}{1 + \alpha}} \right)
   e^{t^{\frac{1 - \alpha}{1 + \alpha}} \tilde{H}_t (y, z)} d y \]
are comparable in size (the rest of the integral is small compared to these
two terms). We can show with similar computations as previously, with
\[ \mathcal{A}_+ (t) \assign \int_{y_+ (z_c, t) - \gamma}^{y_+ (z_c, t) +
   \gamma} e^{t^{\frac{1 - \alpha}{1 + \alpha}} \tilde{H}_t (y, z_c)} d y,
   \mathcal{A}_- (t) \assign \int_{y_- (z_c, t) - \gamma}^{y_- (z_c, t) +
   \gamma} e^{t^{\frac{1 - \alpha}{1 + \alpha}} \tilde{H}_t (y, z_c)} d y, \]
that at leading order we have
\[ t^{\frac{\alpha}{1 + \alpha}} f \left( z t^{\frac{1}{1 + \alpha}}, t
   \right) = \kappa \frac{| y_+ (z_c, t) |^{- \alpha} \mathcal{A}_+ (t) + |
   y_- (z_c, t) |^{- \alpha} \mathcal{A}_- (t)}{\mathcal{A}_+ (t)
   +\mathcal{A}_- (t)} . \]
Although $\tilde{H}_t (y_+ (z_c, t), z_c)$ and $\tilde{H}_t (y_- (z_c, t),
z_c)$ converges to the same constant, it is not clear if
\[ e^{(\tilde{H}_t (y_+ (z_c, t), z_c) - \tilde{H}_t (y_- (z_c, t), z_c))
   t^{\frac{1 - \alpha}{1 + \alpha}}} \]
converges or not. We can only show that for any subsequence of time such that
$t_n^{\frac{\alpha}{1 + \alpha}} f \left( z t_n^{\frac{1}{1 + \alpha}}, t_n
\right)$ converges, the limit must be between $\kappa | y_+^{\ast} (z_c) |^{-
\alpha}$ and $\kappa | y^{\ast}_- (z_c) |^{- \alpha}$.

\subsection{The other asymptotic profiles}\label{ss33}

This subsection is devoted to the proofs of Proposition \ref{3rdcase} to
\ref{IDwhynot}. They follow a similar path than the proof of Theorem
\ref{thfo}. As such, we will only give the main arguments of their proof.

\subsubsection{Proof of Proposition \ref{3rdcase}}\label{ss331}

The key difference here, compared to the proof of Theorem \ref{thfo}, is that
the functions $g$ is different. We define, in this subsection,
\[ g (y) = \left\{\begin{array}{l}
     y - \frac{\kappa}{| y |^{\alpha}} \tmop{if} y > 0\\
     y + \frac{\kappa}{| y |^{\alpha}} \tmop{if} y < 0.
   \end{array}\right. \]
Then, with $g_t (y) = y + t^{\frac{\alpha}{1 + \alpha}} f_0 \left( y
t^{\frac{1}{1 + \alpha}} \right)$, we still have
\[ \| g_t - g \|_{L^{\infty} (\mathbb{R} \backslash [- \varepsilon,
   \varepsilon])} \rightarrow 0 \]
when $t \rightarrow + \infty$. Here, we have $g' (y) > 0$ everywhere on
$\mathbb{R}^{\ast}$, and we thus define
\[ y_+^{\ast} : \mathbb{R} \rightarrow] - \infty, 0 [, y_-^{\ast} : \mathbb{R}
   \rightarrow] 0, + \infty [ \]
as the respective inverses of $g :] - \infty, 0 [\rightarrow \mathbb{R}$ and
$g :] 0, + \infty [\rightarrow \mathbb{R}$.

\

Then, we take $\mu > 0$ small and we decompose in three cases : $z \leqslant
\frac{- 1}{\mu}, z \geqslant \frac{1}{\mu}$ and $z \in \left[ \frac{- 1}{\mu},
\frac{1}{\mu} \right]$. To find where the maximum of $\tilde{H}_t (., z)$ is
reached, the first two cases can be treated like the case $z \geqslant
\frac{1}{\mu}$ in the proof of Theorem \ref{thfo}, where the maximum will be
reached respectively at $y_+ (z, t)$ and $y_- (z, t)$, and in the middle this
can be treated as the case $g (y_0) + \mu \leqslant z \leqslant \frac{1}{\mu}$
in the proof of Theorem \ref{thfo}. We conclude with computations similars as
the ones of subsection \ref{ss325}.

\

Remark that we did not infer any result in the case
\[ f_0 (y) = \frac{\mp \kappa (1 + o_{y \rightarrow \pm \infty} (1))}{| y
   |^{\alpha}}, \]
and this is because our approach does not work here. The difficulty there is
that, with the function
\[ g (y) = \left\{\begin{array}{l}
     y + \frac{\kappa}{| y |^{\alpha}} \tmop{if} y > 0\\
     y - \frac{\kappa}{| y |^{\alpha}} \tmop{if} y < 0,
   \end{array}\right. \]
there exists values of $z$ close to $0$ such that the problem $z = g (y)$ does
not admit any solution. This means that the solutions of $z = g_t (y)$ are
only for values of $y$ close to $0$ when $t \rightarrow + \infty$ in that
case.

As such, we can show the convergence of $t^{\frac{\alpha}{1 + \alpha}} f
\left( z t^{\frac{1}{1 + \alpha}}, t \right)$ for large $| z |$, but not in
the middle, where it is unclear if the profile converge, explode or oscillate.

\subsubsection{Proof of Proposition \ref{lng}}

We consider here the problem
\[ \left\{\begin{array}{l}
     \partial_t f - \partial_x^2 f + f \partial_x f = 0\\
     f_0 (x) = \frac{\kappa}{| y |^{\alpha} \ln^{\beta} (| y |)} (1 + o_{y
     \rightarrow \pm \infty} (1)) .
   \end{array}\right. \]
We write
\[ \lambda (t) f (\mu (t) z, t) = \frac{\int_{\mathbb{R}} \lambda (t) f_0 (\mu
   (t) z) e^{\frac{\mu^2 (t)}{t} \tilde{H}_t (y, z)} d y}{\int_{\mathbb{R}}
   \lambda e^{\frac{\mu^2 (t)}{t} \tilde{H}_t (y, z)} d y} \]
where $\mu (t), \lambda (t) > 0$ are functions going to $+ \infty$ when $t
\rightarrow + \infty$, and
\[ \tilde{H}_t (y, z) = - \frac{(z - y)^2}{4} - \frac{1}{2} \frac{t}{\mu^2
   (t)} \int_0^{y \mu (t)} f_0 (u) d u. \]
We have
\[ \partial_y \tilde{H}_t (y, z) = \frac{1}{2} \left( z - y - \frac{t}{\mu
   (t)} f_0 (y \mu (t)) \right) . \]
For a fixed $y > 0$, we have
\[ \frac{t}{\mu (t)} f_0 (y \mu (t)) \sim \frac{\kappa t}{\mu (t)^{1 + \alpha}
   | y |^{\alpha} \ln^{\beta} (y \mu (t))} \sim \frac{\kappa t}{\mu (t)^{1 +
   \alpha} \ln^{\beta} (\mu (t)) | y |^{\alpha}} \]
when $t \rightarrow + \infty$. We therefore take $\mu (t)$ such that
\[ \mu (t)^{1 + \alpha} \ln^{\beta} (\mu (t)) = t \]
and $\lambda (t) = \frac{t}{\mu (t)}$. Remark that then,
\[ \mu (t) \sim t^{\frac{1}{1 + \alpha}} \left( \frac{1 + \alpha}{\ln (t)}
   \right)^{\frac{\beta}{1 + \alpha}} \]
when $t \rightarrow + \infty$.

\

Now, we define as in the other cases
\[ g (y) = y + \frac{\kappa}{| y |^{\alpha}} \]
and
\[ g_t (y) = y + \frac{t}{\mu (t)} f_0 (y \mu (t)) . \]
We check that for any $\varepsilon > 0$,
\[ \| g_t - g \|_{C^1 \left( \left[ - \frac{1}{\varepsilon},
   \frac{1}{\varepsilon} \right] \backslash [- \varepsilon, \varepsilon]
   \right)} \rightarrow 0. \]
when $t \rightarrow + \infty$. Remark that here the convergence is uniform on
$\left[ - \frac{1}{\varepsilon}, \frac{1}{\varepsilon} \right] \backslash [-
\varepsilon, \varepsilon]$ and not on $\mathbb{R} \backslash [- \varepsilon,
\varepsilon]$ as in the other cases, because for $| y | \geqslant
\varepsilon$,
\[ f_0 (y \mu (t)) \sim \frac{\kappa}{| y |^{\alpha}} \mu (t)^{- \alpha} (\ln
   (| y |) + \ln (\mu (t)))^{- \beta} \sim \frac{\kappa}{| y |^{\alpha}} \mu
   (t)^{- \alpha} \ln^{- \beta} (\mu (t)) \]
when $t \rightarrow + \infty$, but not uniformly in $| y |$ large.

\

We define $y^{\ast}_{\pm}$ to be the same branch as in subsection
\ref{sec41}. The steps of the proof are the same, the only difference is that
we are working on $\left[ - \frac{1}{\varepsilon}, \frac{1}{\varepsilon}
\right]$ from the start instead than on $\mathbb{R}$.

\subsubsection{Proof of Proposition \ref{IDwhynot}}

We give more details for the proof of this proposition, because it diverges at
some point to the proof of Theorem \ref{thfo}. This is because for some values
of $z$, the maximum of $\partial_y \tilde{H}_t (y, z)$ will be reached for $y$
close to $0$, where we have to be careful about the convergences.

\

We recall that
\[ t^{\frac{\alpha}{1 + \alpha}} f \left( z t^{\frac{1}{1 + \alpha}}, t
   \right) = \frac{\int_{\mathbb{R}} t^{\frac{\alpha}{1 + \alpha}} f_0 \left(
   y t^{\frac{1}{1 + \alpha}} \right) e^{t^{\frac{1 - \alpha}{1 + \alpha}}
   \tilde{H}_t (y, z)} d y}{\int_{\mathbb{R}} e^{t^{\frac{1 - \alpha}{1 +
   \alpha}} \tilde{H}_t (y, z)} d y}, \]
\[ \tilde{H}_t (y, z) = - \frac{(z - y)^2}{4} - \frac{1}{2} t^{- \frac{1 -
   \alpha}{1 + \alpha}} \int_0^{y t^{\frac{1}{1 + \alpha}}} f_0 (u) d u, \]
\[ \partial_y \tilde{H}_t (y, z) = \frac{1}{2} \left( z - y -
   t^{\frac{\alpha}{1 + \alpha}} f_0 \left( y t^{\frac{1}{1 + \alpha}} \right)
   \right) \]
and
\[ \partial_y^2 \tilde{H}_t (y, z) = \frac{1}{2} \left( - 1 - t f_0' \left( y
   t^{\frac{1}{1 + \alpha}} \right) \right) . \]
we define
\[ g (y) = \left\{\begin{array}{l}
     y + \frac{\kappa}{| y |^{\alpha}} \tmop{if} y > 0\\
     y \tmop{if} y < 0
   \end{array}\right. \]
and
\[ g_t (y) = y + t^{\frac{\alpha}{1 + \alpha}} f_0 \left( y t^{\frac{1}{1 +
   \alpha}} \right) . \]
For any $\varepsilon > 0$, as in the previous cases, we have
\[ \| g_t - g \|_{C^1 (\mathbb{R} \backslash [- \varepsilon, \varepsilon])}
   \rightarrow 0 \]
when $t \rightarrow + \infty$. We define $y_0 = (\kappa \alpha)^{\frac{1}{1 +
\alpha}} > 0$, the only point where $g' (y_0) = 0$ on $\mathbb{R}^{+ \ast}$.
On $] y_0, + \infty [$ we have $g' > 0$, hence we define $y_+^{\ast} :] g
(y_0), + \infty [\rightarrow] y_0, + \infty [$ as the inverse of $g :] y_0, +
\infty [\rightarrow] g (y_0), + \infty [$. It is the exact same definition as
in subsection \ref{sec41}.

For $y < 0$, since $g (y) = y$, its inverse is simply the identity there.

\

Let us start with the case $z < 0$. Remark that if $y \geqslant 0$, then by
the hypotheses of Proposition \ref{IDwhynot} we have $g (y) \geqslant 0$.
Furthermore, on $] - \infty, 0 [$, we have
\[ g'_t (y) = 1 + t f_0' \left( y t^{\frac{1}{1 + \alpha}} \right) \geqslant
   1, \]
and since $g_t (0) = t^{\frac{\alpha}{1 + \alpha}} f_0 (0) \rightarrow +
\infty$ when $t \rightarrow + \infty$, $g_t (y) \rightarrow - \infty$ when $y
\rightarrow - \infty$, the problem $z = g_t (y)$ admits exactly one solution
for $t$ large enough (depending on $z$). The maximum, denoted $\mathbbmss{y}_{t,
z} < 0$. Given $\nu > 0$ small enough, this branch continues up to $z =
\frac{1}{\nu}$ if $t$ is large enough (depending on $\nu$). For $z < 0$, it
satisfies $\mathbbmss{y}_{t, z} \rightarrow z$ when $t \rightarrow + \infty$.

\

Indeed, for $y < 0$ we have
\[ g_t (y) \sim y + \kappa t^{\frac{\alpha - \beta}{1 + \alpha}} | y |^{-
   \beta} \rightarrow y \]
when $t \rightarrow + \infty$ since $\beta > \alpha$. Now,
\[ \tilde{H}_t (\mathbbmss{y}_{t, z}, z) \rightarrow 0 \]
when $t \rightarrow + \infty$. Therefore there exists $\nu (z) > 0$ such that,
for $y \in \mathbb{R} \backslash \left[ z - \frac{| z |}{2}, z + \frac{| z
|}{2} \right]$ we have
\[ \tilde{H}_t (y, z) < - \nu (z) . \]
Since for $y \in \left[ z - \frac{| z |}{2}, z + \frac{| z |}{2} \right]$ we
have
\[ t^{\frac{\alpha}{1 + \alpha}} f_0 \left( y t^{\frac{1}{1 + \alpha}} \right)
   \leqslant K (z) t^{\frac{\alpha - \beta}{1 + \alpha}} \rightarrow 0 \]
when $t \rightarrow + \infty$, we deduce, with arguments similar as in
subsection \ref{ss325}, that
\[ t^{\frac{\alpha}{1 + \alpha}} f \left( z t^{\frac{1}{1 + \alpha}}, t
   \right) \rightarrow 0 \]
when $t \rightarrow + \infty$, uniformly on $] - \infty, - \varepsilon [$ for
any $\varepsilon > 0$.

\

We now focus on the case $0 < z < g (y_0) + \nu$ for some small $\nu > 0$.
There, the problem $g_t (y) = z$ admits exactly one solution in $] - \infty, 0
[$ (it might have others in $[0, + \infty [$), that we still denote
$\mathbbmss{y}_{t, z} < 0$ (as it is the continuation of it). Here, we have that
$t^{\frac{1}{1 + \alpha}} \mathbbmss{y}_{t, z} \rightarrow - \infty$ when $t
\rightarrow + \infty$, since otherwise, $\mathbbmss{y}_{t, z}$ is bounded and
\[ 0 < z -\mathbbmss{y}_{t, z} = t^{\frac{\alpha}{1 + \alpha}} f_0 \left(
   \mathbbmss{y}_{t, z} t^{\frac{1}{1 + \alpha}} \right) \rightarrow + \infty \]
when $t \rightarrow + \infty$, which is a contradiction. This leads to
\[ \mathbbmss{y}_{t, z} t^{- \frac{\alpha - \beta}{(1 + \alpha) \beta}}
   \rightarrow - \left( \frac{z}{\kappa} \right)^{\frac{1}{\beta}} \]
when $t \rightarrow + \infty$. In particular $\mathbbmss{y}_{t, z} \rightarrow 0$
when $t \rightarrow + \infty$.

Now, for $0 < z < g (y_0) + \nu$, the maximum is reached either at
$\mathbbmss{y}_{t, z}$, or near $A$. It cannot be reached near $A$ because of
arguments similar as in subsection \ref{ss324}.

Remark that
\[ t^{\frac{\alpha}{1 + \alpha}} f_0 \left( \mathbbmss{y}_{t, z} t^{\frac{1}{1 +
   \alpha}} \right) = z -\mathbbmss{y}_{t, z} \rightarrow z \]
when $t \rightarrow + \infty$, and for $y \in [\mathbbmss{y}_{t, z} - t^{-
\lambda}, \mathbbmss{y}_{t, z} + t^{- \lambda}], \lambda > 0$, we have
\begin{eqnarray*}
  &  & \left| t^{\frac{\alpha}{1 + \alpha}} f_0 \left( y t^{\frac{1}{1 +
  \alpha}}, t \right) - t^{\frac{\alpha}{1 + \alpha}} f_0 \left(
  \mathbbmss{y}_{t, z} t^{\frac{1}{1 + \alpha}} \right) \right|\\
  & \leqslant & K t^{- \lambda} t \left| f_0' \left( \mathbbmss{y}_{t, z}
  t^{\frac{1}{1 + \alpha}} \right) \right|\\
  & \leqslant & K t^{1 - \lambda - \frac{(1 + \beta) \alpha}{(1 + \alpha)
  \beta}},
\end{eqnarray*}
hence if $\lambda > \frac{\beta - \alpha}{(1 + \alpha) \beta} > 0$ then $1 -
\lambda - \frac{(1 + \beta) \alpha}{(1 + \alpha) \beta} < 0$ and
\[ \left| t^{\frac{\alpha}{1 + \alpha}} f_0 \left( y t^{\frac{1}{1 + \alpha}},
   t \right) - z \right| \leqslant K t^{\frac{\beta - \alpha}{(1 + \alpha)
   \beta} - \lambda} + \left| t^{\frac{\alpha}{1 + \alpha}} f_0 \left(
   \mathbbmss{y}_{t, z} t^{\frac{1}{1 + \alpha}} \right) +\mathbbmss{y}_{t, z}
   \right| \rightarrow 0 \]
when $t \rightarrow + \infty$. We therefore define
\[ \mathcal{D}_{z, t, \varepsilon} \assign \left[ \mathbbmss{y}_{t, z} -
   t^{\frac{\alpha - \beta}{(1 + \alpha) \beta} - \varepsilon}, \mathbbmss{y}_{t,
   z} + t^{\frac{\alpha - \beta}{(1 + \alpha) \beta} - \varepsilon} \right] .
\]
For $y \in \mathcal{D}_{z, t, \varepsilon}$, with arguments similar as in
subsection \ref{ss324}, there exists $C_1, C_2 > 0$ such that
\[ \tilde{H}_t (y, z) \leqslant \tilde{H}_{t, \max} - C_1 t^{\frac{\beta -
   \alpha}{(1 + \alpha) \beta}} (y -\mathbbmss{y}_{t, z})^2 \]
and
\[ \tilde{H}_t (y, z) \geqslant \tilde{H}_{t, \max} - C_2 t^{\frac{\beta -
   \alpha}{(1 + \alpha) \beta}} (y -\mathbbmss{y}_{t, z})^2 . \]
We compute that
\begin{eqnarray*}
  &  & \int_{\mathcal{D}_{z, t, \varepsilon}} e^{t^{\frac{1 - \alpha}{1 +
  \alpha}} \tilde{H}_t (y, z)} d y\\
  & \geqslant & e^{t^{\frac{1 - \alpha}{1 + \alpha}} \tilde{H}_{t, \max}}
  \int_{\mathcal{D}_{z, t, \varepsilon}} e^{- C t^{\frac{1 - \alpha}{1 +
  \alpha} + \frac{\beta - \alpha}{(1 + \alpha) \beta}} (y - \mathbbmss{y}_{t,z})^2} d y\\
  & \geqslant & e^{t^{\frac{1 - \alpha}{1 + \alpha}} \tilde{H}_{t, \max}}
  \int_{- t^{\frac{\alpha - \beta}{(1 + \alpha) \beta} -
  \varepsilon}}^{t^{\frac{\alpha - \beta}{(1 + \alpha) \beta} - \varepsilon}}
  e^{- C t^{\frac{1 - \alpha}{1 + \alpha} + \frac{\beta - \alpha}{(1 + \alpha)
  \beta}} y^2} d y\\
  & \geqslant & t^{\gamma (\alpha, \beta, \varepsilon)} e^{t^{\frac{1 -
  \alpha}{1 + \alpha}} \tilde{H}_{t, \max}}
\end{eqnarray*}
for some function $\gamma (\alpha, \beta, \varepsilon) > 0$ and $t$ large
enough.

\

We continue, for $y \in \mathbb{R} \backslash \mathcal{D}_{z, t,
\varepsilon}$, we have $| y -\mathbbmss{y}_{t, z} | \geqslant t^{\frac{\alpha -
\beta}{(1 + \alpha) \beta} - \varepsilon}$, therefore
\[ \tilde{H}_t (y, z) \leqslant \tilde{H}_{t, \max} - C_1 t^{\frac{\beta -
   \alpha}{(1 + \alpha) \beta}} t^{2 \left( \frac{\alpha - \beta}{(1 + \alpha)
   \beta} - \varepsilon \right)}, \]
leading to (with arguments similar as for the proof of (\ref{123}) to deal
with the integrability)
\begin{eqnarray*}
  &  & \int_{\mathbb{R} \backslash \mathcal{D}_{z, t, \varepsilon}}
  e^{t^{\frac{1 - \alpha}{1 + \alpha}} \tilde{H}_t (y, z)} d y\\
  & \leqslant & t^{\gamma_1 (\alpha, \beta, \varepsilon)} e^{t^{\frac{1 -
  \alpha}{1 + \alpha}} \tilde{H}_{t, \max}} e^{- C t^{\frac{1 - \alpha}{1 +
  \alpha} + \frac{\beta - \alpha}{(1 + \alpha) \beta} + 2 \left( \frac{\alpha
  - \beta}{(1 + \alpha) \beta} - \varepsilon \right)}}
\end{eqnarray*}
for some function $\gamma_1 (\alpha, \beta, \varepsilon) > 0$, and
\[ \frac{1 - \alpha}{1 + \alpha} + \frac{\beta - \alpha}{(1 + \alpha) \beta} +
   2 \left( \frac{\alpha - \beta}{(1 + \alpha) \beta} - \varepsilon \right) =
   \frac{\alpha (1 - \beta)}{(1 + \alpha) \beta} - 2 \varepsilon > 0 \]
if $\beta < 1$ and $\varepsilon > 0$ small enough.

\

We deduce that
\[ \left| t^{\frac{\alpha}{1 + \alpha}} f \left( z t^{\frac{1}{1 + \alpha}}, t
   \right) - \frac{\int_{\mathcal{D}_{z, t, \varepsilon}} t^{\frac{\alpha}{1 +
   \alpha}} f_0 \left( y t^{\frac{1}{1 + \alpha}} \right) e^{t^{\frac{1 -
   \alpha}{1 + \alpha}} \tilde{H}_t (y, z)} d y}{\int_{\mathbb{R}}
   e^{t^{\frac{1 - \alpha}{1 + \alpha}} \tilde{H}_t (y, z)} d y} \right|
   \rightarrow 0 \]
when $t \rightarrow + \infty$, and since on $\mathcal{D}_{z, t, \varepsilon}$
we have $\left| t^{\frac{\alpha}{1 + \alpha}} f_0 \left( y t^{\frac{1}{1 +
\alpha}}, t \right) - z \right| \rightarrow 0$ when $t \rightarrow + \infty$,
We have
\begin{eqnarray*}
  &  & \frac{\int_{\mathcal{D}_{z, t, \varepsilon}} t^{\frac{\alpha}{1 +
  \alpha}} f_0 \left( y t^{\frac{1}{1 + \alpha}} \right) e^{t^{\frac{1 -
  \alpha}{1 + \alpha}} \tilde{H}_t (y, z)} d y}{\int_{\mathbb{R}}
  e^{t^{\frac{1 - \alpha}{1 + \alpha}} \tilde{H}_t (y, z)} d y} - z\\
  & = & \frac{\int_{\mathcal{D}_{z, t, \varepsilon}} \left(
  t^{\frac{\alpha}{1 + \alpha}} f_0 \left( y t^{\frac{1}{1 + \alpha}} \right)
  - z \right) e^{t^{\frac{1 - \alpha}{1 + \alpha}} \tilde{H}_t (y, z)} d
  y}{\int_{\mathbb{R}} e^{t^{\frac{1 - \alpha}{1 + \alpha}} \tilde{H}_t (y,
  z)} d y}\\
  & - & z \frac{\int_{\mathbb{R} \backslash \mathcal{D}_{z, t, \varepsilon}}
  e^{t^{\frac{1 - \alpha}{1 + \alpha}} \tilde{H}_t (y, z)} d
  y}{\int_{\mathbb{R}} e^{t^{\frac{1 - \alpha}{1 + \alpha}} \tilde{H}_t (y,
  z)} d y}\\
  & \rightarrow & 0
\end{eqnarray*}
when $t \rightarrow + \infty$. Here, we have
\[ \tilde{H}_{t, \max} = - \frac{(z -\mathbbmss{y}_{t, z})^2}{4} - \frac{1}{2}
   t^{- \frac{1 - \alpha}{1 + \alpha}} \int_0^{\mathbbmss{y}_{t, z} t^{\frac{1}{1
   + \alpha}}} f_0 (u) d u \]
and since $\mathbbmss{y}_{t, z} t^{- \frac{\alpha - \beta}{(1 + \alpha) \beta}}
\rightarrow - \left( \frac{z}{\kappa} \right)^{\frac{1}{\beta}}$ when $t
\rightarrow + \infty$, we have
\[ \tilde{H}_{t, \max} \rightarrow - \frac{z^2}{4} \]
when $t \rightarrow + \infty$.

\

Let us now look at the case $\frac{1}{\nu} > z > g (y_0) + \nu$. There, there
are three solutions of $z = g_t (y)$ for $t$ large enough (depending on
$\nu$). One is the continuation of the branch $\mathbbmss{y}_{t, z}$, and the
other two are $y_+ (z, t)$ and $y_m (z, t)$, defined as in the proof of
Theorem \ref{thfo}. We check, as previously, that $y_m$ is a local minimizer
of $\tilde{H}_t (., z)$, so the maximum is either reached at $\mathbbmss{y}_{t,
z}$ or $y_+ (z, t)$.

Defining here
\[ \mathfrak{h}_+ (z) = \tilde{H}_t (y_+ (z, t), z) \text{ and } \mathfrak{h}_-
   (z) = \tilde{H}_t (\mathbbmss{y}_{t, z}, z), \]
we have
\[ \mathfrak{h}_+' (z) = \frac{- 1}{2} (z - y_+ (z, t)), \mathfrak{h}_-' (z) =
   \frac{- 1}{2} (z -\mathbbmss{y}_{t, z}), \]
hence
\[ \mathfrak{h}_+' (z) -\mathfrak{h}_-' (z) = \frac{y_+ (z, t) -\mathbbmss{y}_{t,
   z}}{2} > 0 \]
if $t$ is large enough since $y_+ (z, t) > C (\nu)$ and $\mathbbmss{y}_{t, z}
\rightarrow 0$ when $t \rightarrow + \infty$. Therefore, following the proof
of Theorem \ref{thfo}, there exists exactly one point $z_e (t) > g (y_0)$,
converging when $t \rightarrow + \infty$ such that the maximum is reached at
$\mathbbmss{y}_{t, z}$ for $z < z_e (t)$, and is reached at $y_+ (z, t)$ for $z >
z_e (t)$.

\

Finally, the last case, namely $z > \frac{1}{\nu}$, can be treated as in the
proof of Theorem \ref{thfo}.

\appendix\section{Proof of Propositions \ref{prop12} and
\ref{prop15}}\label{sec3}

\begin{proof}
  We check easily these two propositions if $t \leqslant 1$. We suppose now
  that $t \geqslant 1$. We recall that the solution of the heat equation is
  \[ f (x, t) = \frac{1}{\sqrt{4 \pi t}} \int_{\mathbb{R}} f_0 (y) e^{-
     \frac{(x - y)^2}{4 t}} d y. \]
  We supposed here that $\frac{\kappa_2}{(1 + | y |)^{\alpha}} \leqslant f_0
  (y) \leqslant \frac{\kappa_1}{(1 + | y |)^{\alpha}}$. We estimate
  \[ \frac{1}{\sqrt{4 \pi t}} \int_{| y | \leqslant \sqrt{t}} f_0 (y) e^{-
     \frac{(x - y)^2}{4 t}} d y \leqslant \frac{K}{\sqrt{1 + t}} \int_{| y |
     \leqslant \sqrt{t}} \frac{1}{(1 + | y |)^{\alpha}} \leqslant \frac{K}{(1
     + t)^{\alpha / 2}} \]
  and
  \begin{eqnarray*}
    \frac{1}{\sqrt{4 \pi t}} \int_{| y | \geqslant \sqrt{t}} f_0 (y) e^{-
    \frac{(x - y)^2}{4 t}} d y & \leqslant & \frac{K}{\sqrt{t}} \sup_{| y |
    \geqslant \sqrt{t}} (f_0) \int_{\mathbb{R}} e^{- \frac{(x - y)^2}{4 t}} d
    y\\
    & \leqslant & \frac{K}{\sqrt{t}} (1 + t)^{- \alpha / 2} \sqrt{t}\\
    & \leqslant & \frac{K}{(1 + t)^{\frac{\alpha}{2}}},
  \end{eqnarray*}
  leading to the upper bound on $\| f (., t) \|_{L^{\infty}}$. Finally,
  \[ f_{} (0, t) \geqslant \frac{K}{\sqrt{4 \pi t}} \int_{| y | \leqslant
     \sqrt{t}} \frac{e^{- \frac{y^2}{4 t}}}{(1 + | y |)^{\alpha}} d y
     \geqslant \frac{K}{\sqrt{4 \pi t}} \int_{| y | \leqslant \sqrt{t}}
     \frac{1}{(1 + | y |)^{\alpha}} \geqslant \frac{K}{(1 + t)^{\alpha / 2}} .
  \]
  This completes the proof of Proposition \ref{prop12}. Now, if we define
  \[ G (x, t) = \frac{1}{\sqrt{4 \pi t}} e^{- \frac{x^2}{4 t}}, \]
  by standard computations, we have that for any $n, k \in \mathbb{N}$, there
  exists $K_{n, k} > 0$ such that
  \[ | \partial_t^n \partial_x^k G (x, t) | \leqslant \frac{K_{k,
     n}}{t^{\frac{1}{2} + n + \frac{k}{2}}} e^{- \frac{x^2}{8 t}} . \]
  Therefore,
  \[ | \partial_t^n \partial_x^k f (x, t) | \leqslant \frac{K_{k,
     n}}{t^{\frac{1}{2} + n + \frac{k}{2}}} \int_{\mathbb{R}} f_0 (y) e^{-
     \frac{(x - y)^2}{8 t}} d y \]
  and Proposition \ref{prop15} follows from the same upper bound estimate as
  previously.
\end{proof}

\end{document}